\documentclass[aop,seceqn,MSNbibl,citesort,dvips]{arximspdf}

\doi{10.1214/09-AOP522}
\volume{38}
\issue{4}
\pubyear{2010}
\firstpage{1444}
\lastpage{1491}

\makeatletter

\newcommand{\iint}{\int\!\!\int}

\renewcommand{\mathit}{\mathrm}

\newtheorem{theorem}{Theorem}[section]
\newtheorem{lem}{Lemma}[section]
\newtheorem{cor}[lem]{Corollary}
\newtheorem{prop}[theorem]{Proposition}
\newproclaim{defn}[theorem]{Definition}

\newproclaim{rem}[theorem]{Remark}
\newtheorem{remark}{Remark}[section]

\makeatother

\begin{document}
\begin{frontmatter}

\title{The circular law for random matrices}
\runtitle{The circular law for random matrices}

\begin{aug}
\author[A]{\fnms{Friedrich} \snm{G\"{o}tze}\corref{}\thanksref{t2}\ead[label=e1]{goetze@math.uni-bielefeld.de}} and
\author[B]{\fnms{Alexander} \snm{Tikhomirov}\thanksref{t2,t3}\ead[label=e2]{tichomir@math.uni-bielefeld.de}}
\runauthor{F. G\"{o}tze and A. Tikhomirov}
\affiliation{University of Bielefeld and St. Petersburg State University}
\address[A]{Faculty of Mathematics\\
University of Bielefeld\\
Bielefeld\\
Germany\\
\printead{e1}}
\address[B]{Faculty of Mathematics and Mechanics\\
Sankt-Peterburg State University\\
St. Petersburg\\
Russia\\
\printead{e2}}
\end{aug}

\thankstext{t2}{Supported in part by DFG Project G0-420/5-1 and CRC 701.}
\thankstext{t3}{Supported in part by RFBF Grant N 09-01-12180-ofi\_m,
by RFBR--DFG Grant N 09-01-91331-NNIO\_a.}

\pdfauthor{Friedrich Gotze, Alexander Tikhomirov}

\received{\smonth{11} \syear{2007}}
\revised{\smonth{12} \syear{2008}}

%
\begin{abstract}
We consider the joint distribution of real and imaginary parts of
eigenvalues of random matrices with independent entries with mean
zero and unit variance. We prove the convergence of this distribution to
the uniform distribution on the unit disc without assumptions on the
existence of a density for the distribution of entries. We assume
that the entries have a finite moment of order larger than two and consider
the case of sparse matrices.

The results are based on previous work of Bai, Rudelson and the authors
extending those results to a larger class of sparse matrices.
\end{abstract}

%

\end{frontmatter}

\section{Introduction}

Let $X_{jk}, 1\le j, k<\infty$,
be complex random variables
with $\mathbf{E}X_{jk}=0$ and
$\mathbf{E}|X_{jk}|^2=1$.
For a fixed $n\geq1$,
denote by $\lambda_1,\ldots,\lambda_n$ the eigenvalues of
the $n\times n$ matrix
%
%
\begin{equation}\label{sym}\quad
\mathbf X=(X_n{(j,k)})_{j,k=1}^n,
\qquad X_n{(j,k)}=\frac1{\sqrt n}X_{jk}\qquad\mbox{for }1\leq j,k\leq n,
\end{equation}
and define its empirical spectral distribution function by
%
%
\begin{equation}
G_n(x,y)=\frac1n\sum_{j=1}^nI_{\{\operatorname{Re}\{\lambda_j\}\leq
x, \operatorname{Im}\{\lambda_j\}\le y\}},
\end{equation}
where $I_{\{B\}}$ denotes the indicator of an event $B$. We investigate
the convergence of the expected spectral distribution function
$\mathbf{E}G_n(x,y)$ to the distribution function $G(x,y)$ of the uniform
distribution in the unit disc in $\mathbb R^2$.

%
The main result of our paper is the following:
\begin{theorem}\label{thm0}
Let $\varphi(x)$ denote the function $(\ln(1+|x|))^{19+\eta}$,
$\eta>0$,
arbitrary, small and fixed. Let $X_{jk}, j,k \in
\mathbf N $, denote independent complex random variables with
\[
\mathbf{E}X_{jk}=0,\qquad \mathbf{E}|X_{jk}|^2=1\quad
\mbox{and}\quad \varkappa:=
\sup_{j,k\in\mathbf N} \mathbf{E}
|X_{jk}|^{2}\varphi(X_{jk}) <\infty.
\]
%
Then $\mathbf{E}G_n(x,y)$ converges weakly to the distribution
function $G(x,y)$ as $n\to\infty$.
\end{theorem}

We shall prove the same result for the following class of sparse
matrices. Let
$\varepsilon_{jk}$, $j,k=1,\ldots,n$, denote a triangular array of Bernoulli
random variables (taking values $0,1$ only) which are independent in
aggregate and independent of $(X_{jk})_{j,k=1}^n$ with common success
probability
$p_n:=\Pr\{\varepsilon_{jk}=1\}$ depending on $n$. Consider\vspace*{-1pt} the sequence
of matrices $\mathbf
X^{(\varepsilon)}=\frac1{\sqrt{np_n}}(\varepsilon_{jk}X_{jk})_{j,k=1}^n$.
Let $\lambda_1^{(\varepsilon)},\ldots,\lambda_n^{(\varepsilon)}$
denote the (complex)
eigenvalues of the matrix $\mathbf X^{(\varepsilon)}$ and denote by
$G_n^{(\varepsilon)}(x,y)$ the empirical spectral distribution
function of the matrix $\mathbf X^{(\varepsilon)}$, that is,
%
%
\begin{equation}
G_n^{(\varepsilon)}(x,y):=\frac1n\sum_{j=1}^nI_{\{\operatorname
{Re}\{\lambda_j^{(\varepsilon)}\}\leq x, \operatorname{Im}\{\lambda
_j^{(\varepsilon)}\}\le y\}}.
\end{equation}
%
\begin{theorem}\label{sparse}
For $\eta>0$ define $\varphi(x)=(\ln(1+|x|))^{19+\eta}$. Let
$X_{jk}, j,k \in
\mathbf N $, denote independent complex random variables with
\[
\mathbf{E}X_{jk}=0,\qquad \mathbf{E}|X_{jk}|^2=1\quad
\mbox{and}\quad \varkappa:=
\sup_{j,k\in\mathbf N} \mathbf{E}
|X_{jk}|^{2}\varphi(X_{jk}) <\infty.
\]
%
Assume that there is a $\theta\in(0,1]$ such that $p_n^{-1}=\mathcal
O(n^{1-\theta})$ as $n\to\infty$. 
Then $\mathbf{E}G_n^{(\varepsilon)}(x,y)$ converges weakly to the
distribution function $G(x,y)$ as $n\to\infty$.
\end{theorem}
%
\begin{rem} The crucial problem of the proofs of Theorems \ref{thm0}
and \ref{sparse} is
to bound the smallest singular values $s_n(z)$, respectively,
$s_n^{(\varepsilon)}(z)$ of the shifted matrices $\mathbf X-z\mathbf
I$, respectively, $\mathbf X^{(\varepsilon)}
-z\mathbf I$. (See also \cite{Zei04}, page 1561.) These bounds are
based on the
results obtained by Rudelson and Vershynin in \cite{RV}. In a previous
version of this paper \cite{GT07} we have used the corresponding
results of Rudelson
\cite{rud06} proving the circular law in the case of i.i.d.
sub-Gaussian random variables.
In fact, the results in \cite{GT07} actually imply the circular law
for i.i.d. random variables with $\sup_{j,k}\mathbf
{E}|X_{jk}|^4\le
\varkappa_4<\infty$ in view of the fact (explicitly stated by
Rudelson in
\cite{rud06}) that in his results the sub-Gaussian condition
is needed for the proof of $\Pr\{\|\mathbf X\| > K \}\le C\exp\{-cn\}
$ only.
Restricting oneself to the set $\Omega_n(z)=\{s_n(z)\le c n^{-3};
\|\mathbf
X\| \le K \}$ for the investigation of the smallest singular values, the
inequality $\Pr\{\Omega_n(z)^{c}\}\le c n^{-1/2}$
follows from the results of Rudelson \cite{rud06} \textup{without}
the assumption of sub-Gaussian tails for the matrix $\mathbf X$.
A similar result has been proved by Pan and Zhou in \cite{PanZhou2007}
based on results of Rudelson and Vershynin \cite{RV} and Bai and
Silverstein \cite{Bainn}.

The strong circular law assuming moment condition
of order larger than $2$ only and comparable
sparsity assumptions was proved independently by Tao and Vu in
\cite{TaoVu2007}
based on their results in \cite{TaoVu2007a} in connection with the
multivariate
Littlewood Offord problem.

The approach in this paper though
is based on the fruitful idea of Rudelson and Vershynin to characterize the
vectors leading to small singular values of matrices with independent entries
via ``compressible'' and ``incompressible'' vectors (see~\cite{RV},
Section 3.2, page 15).
For the approximation of the distribution of singular values of
$\mathbf X- z
\mathbf I$
we use a scheme different from the approach used in Bai~\cite{Bai1997}.
\end{rem}

The investigation of the convergence the spectral distribution functions
of real or complex (nonsymmetric and non-Hermitian) random matrices with
independent entries has a long history. Ginibre's \cite{ginibre}, in 1965,
studied the real, complex and quaternion matrices with i.i.d.
Gaussian entries. He derived the joint density for the distribution of
eigenvalues of matrix. Applying Ginibre's formula, Mehta~\cite{Me}, in
1967, determined the density of the expected spectral distribution function
of random matrices with Gaussian entries with independent real and
imaginary parts and deduced the circle law.
Pastur suggested in 1973 the circular law for the general case (see
\cite{Pastur1973}, page 64).
Using the Ginibre results, Edelman \cite{edelman}, in 1997,
proved the circular law for the matrices with i.i.d.
Gaussian real entries. Rider proved in \cite{Rider2003} and \cite
{Rider2006} results about the spectral radius and about linear
statistics of
eigenvalues of non-Hermitian matrices with Gaussian entries.

Girko \cite{Girko1984a}, in 1984, investigated the circular law for
general matrices with
independent entries assuming that the distribution of the entries has densities.
As pointed out by Bai
\cite{Bai1997}, Girko's proof had serious gaps. Bai in \cite
{Bai1997} gave a proof of the circular law
for random matrices with independent entries assuming that the entries
had bounded
densities and finite sixth moments.
His result does not cover
the case of the Wigner ensemble and in particular ensembles of matrices
with Rademacher entries. These ensembles are of some interest in
various
applications (see, e.g., \cite{Timm2004}).
Girko's \cite{Girko1984a} approach using families
of spectra of Hermitian matrices for a characterization of the circular law
based on the so-called \textit{V-transform} was fruitful for all later
work. See, for example, Girko's Lemma 1 in \cite{Bai1997}.
In fact, Girko \cite{Girko1984a} was the first who used the logarithmic
potential to prove the circular law.
We shall outline his approach using logarithmic potential theory.
Let $\xi$ denote a random variable uniformly distributed over the unit
disc and independent of the matrix $\mathbf X$. For any $r>0$, consider
the matrix
\[
\mathbf X(r)=\mathbf X-r\xi\mathbf I,
\]
where $\mathbf I$ denotes the identity matrix of order $n$.
Let $\mu_n^{(r)}$ (resp., $\mu_n$) be empirical spectral measure of
matrix $\mathbf X(r)$ (resp., $\mathbf X$)
defined on the complex plane
as empirical measure of the set of eigenvalues of matrix.
We define a
logarithmic potential of the expected spectral measure $\mathbf{E}\mu
_n^{(r)}(ds,dt)$
as
\[
U_{\mu_n}^{(r)}(z)=-\frac1n\mathbf{E}\log\bigl|\det\bigl(\mathbf
X(r)-z\mathbf I\bigr)\bigr|
=-\frac1n\sum{\mathbf{E}\log}|\lambda_j-z-r\xi|,
\]
where $\lambda_1,\ldots,\lambda_n$ are the eigenvalues of the matrix
$\mathbf X$. Note that the expected spectral measure $\mathbf{E}\mu
_n^{(r)}$ is the
convolution of the measure $\mathbf{E}\mu_n$ and the uniform
distribution on the disc
of radius $r$ (see Lemma \ref{ap1} in the \hyperref[app]{Appendix} for details).
\begin{lem}\label{sma1} Assume that the sequence $\mathbf{E}\mu
_n^{(r)}$ converges weakly to a measure $\mu$
as $n\to\infty$ and $r\to0$.
Then
\[
\mu=\lim_{n\to\infty}\mathbf{E}\mu_n.
\]
\end{lem}
\begin{pf}
Let $J$ be a random variable which is uniformly distributed on the set
$\{1,\ldots,n\}$ and independent of the matrix $\mathbf X$. We may
represent the measure $\mathbf{E}\mu_n^{(r)}$ as the distribution of
a random variable
$\lambda_J+r\xi$ where $\lambda_J$ and $\xi$ are independent. Computing
the characteristic function of this measure and passing
first to the limit with respect to $n\to\infty$ and then with respect to
$r\to0$ (see also Lemma \ref{ap2} in the \hyperref[app]{Appendix}), we
conclude the result.
\end{pf}

Now we may fix
$r>0$ and consider the measures $\mathbf{E}\mu_n^{(r)}$. They have bounded
densities. Assume that the measures $\mathbf{E}\mu_n$ have
supports in a fixed compact set and that $\mathbf{E}\mu_n$ converges
weakly to a
measure $\mu$. Applying Theorem 6.9 (Lower envelope theorem) from
\cite{saff}, page 73 (see also Section 3.8 in the
Appendix), we obtain that under these assumptions
\[
\liminf_{n\to\infty}U_{\mu_n}^{(r)}(z)=U^{(r)}(z),
\]
quasi-everywhere in $\mathbb C$ (for the definition of
``\textit{quasi-everywhere}'' see, e.g., \cite{saff}, page 24).
Here $U^{(r)}(z)$ denotes the logarithmic potential of the measure
$\mu^{(r)}$ which is the convolution of a measure
$\mu$ and of the uniform distribution on the disc of radius $r$.
Furthermore, note that $U^{(r)}(z)$
may be represented as
\[
U^{(r)}(z_0)=\frac2{ r^2}\int_0^rvL(\mu;z_0,v)\,dv,
\]
where
%
%
\begin{equation}
L(\mu;z_0,v)=\frac1{2\pi}\int_{-\pi}^{\pi}U_{\mu}(z_0+v\exp\{
i\theta\})\,d\theta
\end{equation}
and
%
%
\begin{equation}
U_{\mu}(z)={\int\ln}|\zeta-z|\,d\mu(\zeta).
\end{equation}

Applying Theorem 1.2 in \cite{saff}, page 84,
we get
\[
\lim_{r\to0}U_{\mu}^{(r)}(z)=U_{\mu}(z).
\]
Let $s_1(\mathbf X)\ge\cdots\ge s_n(\mathbf X)$ denote the singular
values of the
matrix $\mathbf X$.\vspace*{2pt}

Since $\mathbf{E}\frac{1}{n}\operatorname{Tr}\mathbf X\mathbf X^*=1$ the
sequence of
measures $\mathbf{E}\mu_n$ is weakly relatively compact. These
results imply that for any $\eta>0$ we may restrict the
measures $\mathbf{E}\mu_n$
to some compact set $K_{\eta}$ such that $\sup_n\mathbf{E}\mu
_n(K_{\eta}^{(c)})<\eta$. Moreover, Lemma \ref{compact} implies the
existence of a compact $K$ such that\break $\lim_{n\to\infty}\sup
_n\mathbf{E}\mu_n(K^{(c)})=0$. If we take some subsequence of the
sequence of restricted measures $\mathbf{E}\mu_n$ which converges to
some measure $\mu$,
then\break $\liminf_{n\to\infty}U_{\mu_n}^{(r)}(z)=U_{\mu
}^{(r)}(z)$, $r>0$, and
$\lim_{r\to0}U_{\mu}^{(r)}(z)=U_{\mu}(z)$. If we prove that
$\liminf_{n\to\infty}U_{\mu_n}^{(r)}(z)$ exists
and $U_{\mu}(z)$ is equal to the logarithmic potential
corresponding the uniform distribution on the unit disc
[see Section \ref{property}, equality (\ref{main})], then the
sequence of
measures $\mathbf{E}\mu_n$ weakly converges to the uniform
distribution on the
unit disc. Moreover, it is enough to prove that for some sequence
$r=r(n)\to0$, $\lim_{n\to\infty}U_{\mu_n}^{(r)}(z)=U_{\mu}(z)$.

Furthermore, let $s_1^{(\varepsilon)}(z,r)\ge\cdots\ge
s_n^{(\varepsilon)}(z,r)$
denote the singular values of matrix $\mathbf X^{(\varepsilon
)}(z,r)=\mathbf
X^{(\varepsilon)}(r)-z\mathbf I$. We shall investigate the logarithmic
potential $U_{\mu_n}^{(r)}(z)$.
Using elementary properties of singular values (see, e.g.,
\cite{Goh69}, Lemma 3.3, page~35), we may represent the function
$U_{\mu_n}^{(r)}(z)$ as follows:
\[
U_{\mu_n}^{(r)}(z)=-\frac1n\sum_{j=1}^n\mathbf{E}\log
{s_j^{(\varepsilon)}(z,r)}
=-\frac12\int_0^{\infty}\log
x \nu_n^{(\varepsilon)}(dx,z,r),
\]
where $\nu_n^{(\varepsilon)}(\cdot,z,r)$ denotes the expected
spectral measure of the
matrix\break $\mathbf H_n^{(\varepsilon)}(z,r)=(\mathbf
X^{(\varepsilon)}(r)-z\mathbf
I)(\mathbf X^{(\varepsilon)}(r)-z\mathbf I)^*$, which is the
expectation of the counting measure of the set of
eigenvalues of the matrix $\mathbf H_n^{(\varepsilon)}(z,r)$.

In Section \ref{convergence} we investigate
convergence of the measure
$\nu_n^{(\varepsilon)}(\cdot,z):=\break\nu^{(\varepsilon)}(\cdot,z,0)$.
In Section
\ref{property} we study the properties of the limit measures
$\nu(\cdot,z)$. But the crucial problem for the proof of the circular
law is
the so-called ``regularization of the potential.'' We solve this
problem using
bounds for the minimal singular values
of the matrices $\mathbf X^{(\varepsilon)}(z):=\mathbf X^{(\varepsilon)}-z\mathbf I$
based on techniques developed in Rudelson \cite{rud06} and Rudelson
and Vershynin \cite{RV}.
The bounds of minimal singular values of matrices $\mathbf
X^{(\varepsilon)}$
are given in Section \ref{singular} and in the
Appendix,
Theorem \ref{sparse}.
In Section~\ref{proof}
we give the proof of the main theorem. In the \hyperref[app]{Appendix}
we combine
precise statements of relevant results
from potential theory and some auxiliary
inequalities for the resolvent matrices.

In the what follows we shall denote by $C$ and $c$ or $\alpha,\beta,
\delta,\rho, \eta$ (without indices) some general absolute constant
which may be changed from line to line.
To specify a constant we shall use subindices. By $I_A$ we shall denote
the indicator of an event~$A$.
For any matrix $\mathbf G$ we denote the Frobenius norm by $\|\mathbf
G\|_2$, and we denote by $\|\mathbf G\|$ the operator norm.

\section{Convergence of $\nu_n^{(\varepsilon)}(\cdot,z)$}
\label{convergence}

Denote by $F_n^{(\varepsilon)}(x,z)$ the distribution function
of the measure $\nu_n^{(\varepsilon)}(\cdot,z)$, that is,
\[
F_n^{(\varepsilon)}(x,z)=\frac1n\sum_{j=1}^n\mathbf{E}I_{\{
s_j^{(\varepsilon)}{(z)}^2<x\}},
\]
where $s_1^{(\varepsilon)}(z)\ge\cdots\ge s_n^{(\varepsilon)}(z)\ge0$
denote the singular values of the matrix $\mathbf X^{(\varepsilon
)}(z)=\mathbf
X^{(\varepsilon)}-z\mathbf I$. For a
positive random variable $\xi$ and a Rademacher random variable (r.v.)
$\kappa$ consider the transformed r.v. $\widetilde\xi=\kappa\sqrt
{\xi}$.
If $\zeta$ has distribution function ${F}_n^{(\varepsilon)}(x,z)$,
the variable $\widetilde{\zeta}$ has distribution
function $\widetilde{F}_n^{(\varepsilon)}(x,z)$, given by
\[
\widetilde{F}_n^{(\varepsilon)}(x,z)=\tfrac12 \bigl(1+\operatorname{sgn}
\{x\}F_n^{(\varepsilon)}(x^2,z) \bigr)
\]
for all real $x$.
Note that this induces
a one-to-one corresponds between the respective measures
$\nu_n^{(\varepsilon)}(\cdot,z)$ and ${\widetilde{\nu
}}_n^{(\varepsilon)}(\cdot,z)$.
The limit distribution function of ${F}_n^{(\varepsilon)}(x,z)$ as
$n\to\infty$, is denoted
by $F(\cdot,z)$. The corresponding symmetrization $\widetilde{F}(x,z)$
is the limit of $\widetilde{F}_n^{(\varepsilon)}(x,z)$ as $n\to
\infty$. We have
\[
\sup_x\bigl|F_n^{(\varepsilon)}(x,z)-F(x,z)\bigr|=2\sup_x\bigl|\widetilde
F_n^{(\varepsilon)}(x,z)-\widetilde F(x,z)\bigr|.
\]

Denote by $s_n^{(\varepsilon)}(\alpha,z)$ [resp., $s(\alpha,z)$] and
$S_n^{(\varepsilon)}(x,z)$ [resp., $S(x,z)$] the Stieltjes transforms
of the
measures $\nu_n^{(\varepsilon)}(\cdot,z)$ [resp., $\nu(\cdot,z)$] and
$\widetilde{\nu}_n^{(\varepsilon)}(\cdot,z)$ [resp., $\widetilde
{\nu}(\cdot,z)$] correspondingly. Then we have
\[
S_n^{(\varepsilon)}(\alpha, z)=\alpha s_n^{(\varepsilon)}(\alpha^2,
z),\qquad
S(\alpha, z) =\alpha s(\alpha^2, z).
\]
\begin{rem}
As shown in Bai \cite{Bai1997}, the measure
$\nu(\cdot,z)$ has a density $p(x,z)$ with
bounded support. More precisely,
$p(x,z)\le C\max\{1, \frac1{\sqrt x}\}$.
Thus the measure $\widetilde\nu(\cdot,z)$ has bounded support and
bounded density
$\widetilde p(x,z)=|x|p(x^2,z)$.
\end{rem}
%
\begin{theorem}\label{thm05}Let $\mathbf{E}X_{jk}=0$, $\mathbf
{E}|X_{jk}|^2=1$. Assume for some function $\varphi(x)>0$
such that $\varphi(x)\to\infty$ as $x\to\infty$ and such that the
function $x/\varphi(x)$ is nondecreasing we have
%
%
\begin{equation}
\varkappa:=\max_{1\le j,k<\infty}\mathbf{E}|X_{jk}|^{2}\varphi
(X_{jk})<\infty.
\end{equation}
Then
%
%
\begin{equation}
\sup_x\bigl|F_n^{(\varepsilon)}(x,z)-F(x,z)\bigr|\le C\varkappa\bigl(\varphi\bigl(\sqrt
{np_n}\bigr)\bigr)^{-1/6}.
\end{equation}
\end{theorem}
\begin{cor}\label{cor}
Let $\mathbf{E}X_{jk}=0$, $\mathbf{E}|X_{jk}|^2=1$, and
%
%
\begin{equation}
\varkappa=\max_{1\le j,k<\infty}\mathbf{E}|X_{jk}|^{3}<\infty.
\end{equation}
Then
%
%
\begin{equation}
\sup_x\bigl|F_n^{(\varepsilon)}(x,z)-F(x,z)\bigr|\le C({np_n})^{-1/{12}}.
\end{equation}
\end{cor}
%
\begin{pf}
To bound the distance between the distribution functions\break
$\widetilde F_n^{(\varepsilon)}(x,z)$ and $\widetilde F(x,z)$
we investigate the distance between their the Stieltjes transforms.
Introduce the Hermitian $2n\times2n$ matrix
\[
\mathbf W=\pmatrix{\mathbf O_n & \bigl(\mathbf X^{(\varepsilon
)}-z\mathbf
I\bigr)
\cr
\bigl(\mathbf X^{(\varepsilon)}-z\mathbf I\bigr)^* & \mathbf O_n},
\]
where $\mathbf O_n$ denotes $n\times n$ matrix with zero entries.
Using the inverse of the partial matrix (see, e.g., \cite{HoJohn91},
Chapter 08, page 18) it follows that,
for $\alpha=u+iv$, $v>0$,
%
%
\begin{eqnarray}\label{shur}
&&(\mathbf W-\alpha\mathbf I_{2n})^{-1}
=
\left(\matrix{
\alpha\bigl(\mathbf X^{(\varepsilon)}(z)\mathbf
X^{(\varepsilon)}(z)^*-\alpha^2\mathbf I \bigr)^{-1} &
\cr
\mathbf
X^{(\varepsilon)}(z)^* \bigl(\mathbf
X^{(\varepsilon)}(z)\mathbf X^{(\varepsilon)}(z)^*-\alpha^2\mathbf
I \bigr)^{-1} & }\right.
\nonumber\\[-8pt]\\[-8pt]
&&\hspace*{85.2pt}\left.\matrix{
\mathbf
X^{(\varepsilon)}(z) \bigl(\mathbf X^{(\varepsilon)}(z)^*\mathbf
X^{(\varepsilon)}(z)-\alpha^2\mathbf I \bigr)^{-1}
\cr
\alpha\bigl(\mathbf X^{(\varepsilon
)}(z)^*\mathbf
X^{(\varepsilon)}(z)-\alpha^2\mathbf I \bigr)^{-1}}\right),\nonumber
\end{eqnarray}
where $\mathbf X^{(\varepsilon)}(z)=\mathbf X^{(\varepsilon
)}-z\mathbf I$ and
$\mathbf I_{2n}$ denotes the unit matrix of order $2n$. By definition of
$S_n^{(\varepsilon)}(\alpha,z)$, we have
\[
S_n^{(\varepsilon)}(\alpha,z)=\frac1{2n}\mathbf{E}\operatorname
{Tr}(\mathbf W-\alpha\mathbf
I_{2n})^{-1}.
\]
Set $\mathbf R(\alpha,z):=(R_{j,k}(\alpha,z))_{j,k=1}^{2n}
=(\mathbf W-\alpha\mathbf I_{2n})^{-1}$.
It is easy to check that
\[
1+\alpha S_n^{(\varepsilon)}(\alpha,z)=\frac{1}{2n}\mathbf
{E}\operatorname{Tr}\mathbf W\mathbf
R(\alpha,z).
\]
We may rewrite this equality as
%
%
\begin{eqnarray}\label{n1}
&&
1+\alpha S_n^{(\varepsilon)}(\alpha,z)\nonumber\\
&&\qquad=\frac1{2n\sqrt
{np_n}}\sum_{j,k=1}^n \mathbf{E}
\bigl({\varepsilon_{jk}}X_{jk}R_{k+n,j}(\alpha,z)\nonumber\\[-8pt]\\[-8pt]
&&\hspace*{107.65pt}{}+{\varepsilon
_{jk}}\overline
X_{jk}
R_{k,j+n}(\alpha,z)\bigr)\nonumber\\
&&\qquad\quad{}-\frac{\overline z}{2n}\sum_{j=1}^n\mathbf{E}R_{j,j+n}(\alpha,z)
-\frac{z}{2n}\sum_{j=1}^n\mathbf{E}R_{j+n,j}(\alpha,z).\nonumber
\end{eqnarray}
We introduce the notation
\begin{eqnarray*}
\mathbf A&=&\bigl(\mathbf X^{(\varepsilon)}(z)\mathbf
X^{(\varepsilon)}(z)^*-\alpha^2\mathbf I\bigr)^{-1},\qquad
\mathbf B=\mathbf X^{(\varepsilon)}(z)\mathbf C,\\
\mathbf C&=&\bigl(\mathbf X^{(\varepsilon)}(z)^*\mathbf
X^{(\varepsilon)}(z)-\alpha^2\mathbf I\bigr)^{-1},\qquad \mathbf D=\mathbf
X^{(\varepsilon)}(z)^*\mathbf A.
\end{eqnarray*}
With this notation we rewrite equality (\ref{shur}) as follows:
%
%
\begin{equation}\label{shur1}
\mathbf R(\alpha,z)=(\mathbf W-\alpha\mathbf I_{2n})^{-1}=
\pmatrix{\alpha\mathbf A &
\mathbf B
\cr
\mathbf D &
\alpha\mathbf C}.
\end{equation}
Equalities (\ref{shur1}) and (\ref{n1}) together imply
%
%
\begin{eqnarray}\label{n2}
&&
1+\alpha S_n^{(\varepsilon)}(\alpha,z)\nonumber\\
&&\qquad=\frac1{2n\sqrt{np_n}}\sum_{j,k=1}^n
\mathbf{E}\bigl({\varepsilon_{jk}}X_{jk}R_{k+n,j}(\alpha,z)\nonumber\\[-8pt]\\[-8pt]
&&\hspace*{107.65pt}{}+\varepsilon
_{jk}\overline X_{jk}R_{k,j+n}(\alpha,z)\bigr)\nonumber\\
&&\qquad\quad{}-\frac{z}{2n}\mathbf{E}\operatorname{Tr}\mathbf D-\frac{\overline
z}{2n}\mathbf{E}\operatorname{Tr}\mathbf B.\nonumber
\end{eqnarray}

In what\vspace*{1pt} follows we shall use a simple resolvent equality.
For two matrices
$\mathbf U$ and $\mathbf V$
let $\mathbf R_{U}=(\mathbf U-\alpha\mathbf I)^{-1}$, $\mathbf R_{U+V}=
(\mathbf U+\mathbf V-\alpha\mathbf I)^{-1}$,
then
\[
\mathbf R_{U+V}=\mathbf R_{U}-\mathbf R_{U}\mathbf V\mathbf R_{U+V}.
\]
Let $\{\mathbf e_1,\ldots,\mathbf e_{2n}\}$ denote the canonical
orthonormal basis
in $\mathbb R^{2n}$.
Let $\mathbf W^{(jk)}$ denote the matrix obtained from
$\mathbf W$ by replacing both entries $X_{j,k}$ and
$\overline X_{j,k}$ by 0. In our notation we may write
%
%
\begin{equation}\label{repr1}
\mathbf W=\mathbf W^{(jk)}+\frac1{\sqrt{np_n}}{\varepsilon
_{jk}}X_{jk}\mathbf
e_{j}\mathbf e_{k+n}^T+ \frac1{\sqrt{np_n}}{\varepsilon
_{jk}}\overline
X_{jk}\mathbf e_{k+n}\mathbf e_{j}^T.
\end{equation}
Using this representation and the resolvent equality, we get
%
%
\begin{eqnarray}\label{rezolvent}
\mathbf R&=&\mathbf R^{(j,k)}-\frac1{\sqrt{np_n}}{\varepsilon
_{jk}}X_{jk}\mathbf
R^{(j,k)}\mathbf e_{j} \mathbf e_{k+n}^T\mathbf R\nonumber\\[-8pt]\\[-8pt]
&&{} -\frac1{\sqrt
{np_n}}{\varepsilon_{jk}}\overline X_{jk}\mathbf R^{(j,k)}\mathbf
e_{k+n}\mathbf e_{j}^T \mathbf R.\nonumber
\end{eqnarray}
Here, and in what follows, we omit the arguments $\alpha$ and $z$ in
the notation of resolvent matrices. For any vector $\mathbf a$, let
$\mathbf a^T$
denote the transposed vector $\mathbf a$. Applying the resolvent
equality again,
we obtain
%
%
\begin{eqnarray}
\mathbf R&=&\mathbf R^{(j,k)}-\frac1{\sqrt{np_n}}{\varepsilon
_{jk}}X_{jk}\mathbf
R^{(j,k)}\mathbf e_{j} \mathbf e_{k+n}^T\mathbf R^{(j,k)}
\nonumber\\[-8pt]\\[-8pt]
&&{}-\frac1{\sqrt
{np_n}}{\varepsilon_{jk}}\overline X_{jk}\mathbf R^{(j,k)}\mathbf
e_{k+n}\mathbf e_{j}^T \mathbf R^{(j,k)} +\mathbf T^{(jk)},\nonumber
\end{eqnarray}
where
%
%
\begin{eqnarray}\label{67}
\mathbf T^{(jk)}&=&\frac1{\sqrt{np_n}}{\varepsilon_{jk}}X_{jk}\mathbf
R^{(j,k)}\mathbf e_{j}\mathbf e_{k+n}^T \bigl(\mathbf R^{(j,k)}-\mathbf
R\bigr)\nonumber\\
&&{}+\frac1{\sqrt{np_n}}{\varepsilon_{jk}}X_{jk}\mathbf R^{(j,k)}\mathbf
e_{j}\mathbf e_{k+n}^T
\bigl(\mathbf R^{(j,k)}-\mathbf R\bigr)\nonumber\\[-8pt]\\[-8pt]
&&{}+\frac1{\sqrt{np_n}}{\varepsilon_{jk}}(\overline X_{jk})\mathbf
R^{(j,k)}\mathbf e_{k+n}\mathbf e_{j}^T\bigl(\mathbf R^{(j,k)}-\mathbf
R\bigr)\nonumber\\
&&{}+\frac1{\sqrt{np_n}}{\varepsilon_{jk}}X_{jk}\mathbf R^{(j,k)}\mathbf
e_{k+n}\mathbf e_{j}^T\bigl(\mathbf R^{(j,k)}-\mathbf R\bigr).\nonumber
\end{eqnarray}
This implies
%
%
\begin{eqnarray}\label{66}
\mathbf R_{j,k+n}&=&\mathbf R^{(j,k)}_{j,k+n}-\frac1{\sqrt
{np_n}}{\varepsilon_{jk}}X_{jk} \mathbf R^{(j,k)}_{j,j}\mathbf
R^{(j,k)}_{k+n,k+n} \nonumber\\
&&{}
-\frac1{\sqrt{np_n}}{\varepsilon_{jk}}\overline
X_{jk}\bigl(\mathbf R^{(j,k)}_{j,k+n}\bigr)^2+
\mathbf T^{(j,k)}_{j,k+n},\nonumber\\[-8pt]\\[-8pt]
\mathbf R_{k+n,j}&=&\mathbf R^{(j,k)}_{k+n,j}-\frac1{\sqrt
{np_n}}{\varepsilon_{jk}}X_{jk} \mathbf R^{(j,k)}_{k+n,j}\mathbf
R^{(j,k)}_{j,k+n} \nonumber\\
&&{}
-\frac1{\sqrt{np_n}}{\varepsilon_{jk}}\overline
X_{jk}\mathbf R^{(j,k)}_{k+n,k+n} \mathbf R^{(j,k)}_{j,j}+\mathbf
T^{(j,k)}_{k+n,j}.\nonumber
\end{eqnarray}
Applying this notation to equality (\ref{n2}) and taking into
account that $X_{jk}$ and
$\mathbf R^{(jk)}$ are independent, we get
%
%
\begin{eqnarray}\label{new12}
&&
1+\alpha S_n^{(\varepsilon)}(\alpha,z)+\frac{z}{2n}\operatorname
{Tr}\mathbf D+\frac{\overline z}{2n}
\operatorname{Tr}\mathbf B\nonumber\\
&&\qquad=-\frac1{n^2p_n}\sum_{j,k=1}^n \mathbf{E}
{\varepsilon_{jk}}|X_{jk}|^2R^{(j,k)}_{j,j}
R^{(j,k)}_{k+n,k+n}\nonumber\\[-8pt]\\[-8pt]
&&\qquad\quad{} -
\frac1{n^2p_n}\sum_{j,k=1}^n \mathbf{E}{\varepsilon
_{jk}}\operatorname{Re}(X_{jk}^2)\mathbf{E}
\bigl(R^{(j,k)}_{j,k+n}\bigr)^2
\nonumber\\
&&\qquad\quad{}-\frac1{2n\sqrt{np_n}}\sum_{j,k=1}^n
\mathbf{E}\bigl({\varepsilon_{jk}}X_{jk}T^{(j,k)}_{k+n,j}+{\varepsilon
_{jk}}\overline
X_{jk}T^{(j,k)}_{j,k+n}\bigr).\nonumber
\end{eqnarray}
From (\ref{rezolvent}) it follows immediately that for any
$p,q=1,\ldots,2n$, $j,k=1,\ldots,n$,
%
%
\begin{equation}\label{bound1}
\bigl|R_{p,p}-R^{(j,k)}_{p,p}\bigr|\le\frac{C{\varepsilon_{jk}}|X_{jk}|}
{\sqrt{np_n}}(|R^{jk}_{pj}||R_{k+n,p}|+|R^{jk}_{p,k+n}||R_{jp}|).
\end{equation}
Since ${\sum_{m,l=1}^n} |R_{m,l}|^2\le n/v^2$ and ${\sum_{m,l=1}^n}
|R^{(jk)}_{m,l}|^2\le n/v^2$, equality (\ref{66}) implies
%
%
\begin{equation}\label{bound1a}
\frac1{n^2}\sum_{j,k=1}^n\mathbf{E}\bigl|R^{(j,k)}_{j,k+n}\bigr|^2\le
\frac{C}{nv^4}.
\end{equation}
By definition (\ref{67}) of $\mathbf T^{(j,k)}$, applying standard
resolvent properties,
we obtain the following bounds,
for any $z=u+iv, v>0$,
%
%
\begin{equation}\label{bound10}
\frac1{n\sqrt
{np_n}}\sum_{j,k=1}^n\mathbf{E}{\varepsilon
_{jk}}|X_{jk}|\bigl|T^{(j,k)}_{j,k+n}\bigr|\le
\frac{C\varkappa}{v^3\varphi(\sqrt{np_n})}.
\end{equation}
For the proof of this inequality see Lemma \ref{ap} in the \hyperref
[app]{Appendix}.
Using the last inequalities we obtain, that for $v>0$
%
%
\begin{eqnarray}\label{41}\qquad
&&\Biggl|\frac1n\sum_{j=1}^n\mathbf{E}
R_{jj}\frac1n\sum_{k=1}^nR_{k+n,k+n}
-\frac1{n^2}
\sum_{j=1}^n\sum_{k=1}^n\mathbf{E}R^{(jk)}_{jj}
R^{(jk)}_{k+n,k+n} \Biggr|\nonumber\\
&&\qquad\le\frac{C}{n^2\sqrt
{np_n}v}\sum_{j=1}^n \sum_{k=1}^n \mathbf{E}\varepsilon_{jk}|X_{jk}|
\bigl(\bigl|R^{(jk)}_{jj}\bigr||R_{k+n,j}| +
\bigl|R^{(jk)}_{j,k+n}\bigr||R_{jj}|\bigr)\\
&&\qquad \le\frac{C}{nv^4}.\nonumber
\end{eqnarray}

Since $\frac1n\sum_{j=1}^nR_{jj}=\frac1n\sum_{k=1}^nR_{k+n,k+n}
=\frac1{2n}\operatorname{Tr}\mathbf R(\alpha,z)$,
we obtain
%
%
\begin{equation}\label{prev}
\Biggl|\frac1{n^2}\sum_{j=1}^n\sum_{k=1}^n\mathbf{E}R^{(jk)}_{jj}
R^{(jk)}_{k+n,k+n}-\mathbf{E}\biggl(\frac1{2n}\operatorname{Tr}\mathbf
R(\alpha,z)\biggr)^2 \Biggr|\le\frac
{C}{nv^4}.
\end{equation}
Note that for any Hermitian random matrix $\mathbf W$ with independent
entries on and above the
diagonal we have
%
%
\begin{equation}\label{variation}
\mathbf{E} \biggl|\frac1n\operatorname{Tr}\mathbf R(\alpha
,z)-\mathbf{E}\frac1n\operatorname{Tr}\mathbf R(\alpha,z)
\biggr|^2\le
\frac{C}{nv^2}.
\end{equation}
The proof of this inequality is easy and due to a martingale-type
expansion already used by Girko. Inequalities (\ref{prev}) and (\ref
{variation}) together imply
that for $v>0$
%
%
\begin{equation}\label{ne6}
\Biggl|\frac1{n^2}\sum_{j=1}^n\sum_{k=1}^n\mathbf{E}R^{(jk)}_{jj}
R^{(jk)}_{k+n,k+n}-\bigl(S_n^{(\varepsilon)}(\alpha,z)\bigr)^2 \Biggr|\le\frac
{C}{ nv^4}.
\end{equation}
Denote by $r(\alpha,z)$ some generic function with $|r(\alpha,z)|\le
1$ which may vary from line to line.
We may now rewrite equality (\ref{n2}) as follows:
%
%
\begin{eqnarray}\label{n3}
&&1+\alpha S_n^{(\varepsilon)}(\alpha,z)+\bigl(S_n^{(\varepsilon)}(\alpha
,z)\bigr)^2\nonumber\\[-8pt]\\[-8pt]
&&\qquad= -\frac{z}{2n}\mathbf{E}\operatorname{Tr}\mathbf
D-\frac{\overline z}{2n}\mathbf{E}\operatorname{Tr}\mathbf B+ \frac
{r(\alpha,z)}{ v^3\varphi(\sqrt{np_n})},\nonumber
\end{eqnarray}
where
$v>c\varphi(\sqrt{np_n})/n$.

We now investigate the functions
$T(\alpha,z)=\frac1n\mathbf{E}\operatorname{Tr}\mathbf B$ and
$V(\alpha,z)=\frac1n\mathbf{E}\operatorname{Tr}\mathbf D$.
Since the arguments for both functions are similar we provide it for the
first one only. By definition of the matrix $\mathbf B$, we have
\[
\operatorname{Tr}\mathbf B=\frac1{\sqrt
{np_n}}\sum_{j,k=1}^n\varepsilon_{jk}X_{j,k}\bigl(\mathbf X^{(\varepsilon
)}(z)^*\mathbf
X^{(\varepsilon)}(z) -\alpha^2\mathbf I\bigr)^{-1}_{kj}-z\operatorname
{Tr}\mathbf C.
\]
According to equality (\ref{shur1}), we have
\[
\operatorname{Tr}\mathbf B=\frac1{\alpha\sqrt
{np_n}}\sum_{j,k=1}^n\varepsilon
_{jk}X_{j,k}R_{k+n,j+n}-z\operatorname{Tr}\mathbf C.
\]

Using the resolvent equality (\ref{rezolvent}) and Lemma \ref{ap}, we get,
for $v>c\times\break\varphi(\sqrt{np_n})/n$
%
%
\begin{equation}\label{n7}\qquad\quad
T(\alpha,z)=-\frac1{\alpha n^2}\sum_{j,k=1}^n\mathbf{E}R^{(jk)}_{k+n,k+n}
R^{(jk)}_{j,j+n}-\frac{z}{\alpha}S_n^{(\varepsilon)}(\alpha,z)
+\frac{C\varkappa r(\alpha, z)}{v^3\varphi(\sqrt{np_n})}.
\end{equation}

Similar to (\ref{ne6}) we obtain
%
%
\begin{equation}\label{n8}
\Biggl|\frac1{n^2}\sum_{j,k=1}^n\mathbf{E}R^{(jk)}_{j,j+n}R^{(jk)}_{k+n,k+n}-
T(\alpha,z)S_n^{(\varepsilon)}(\alpha,z)\Biggr|\le\frac C{nv^4}.
\end{equation}
Inequalities (\ref{n7}) and (\ref{n8}) together imply, for
$v>c\varphi(\sqrt{np_n})/n$,
%
%
\begin{equation}\label{ne8}
T(\alpha,z)=-\frac{zS_n^{(\varepsilon)}(\alpha,z)}{\alpha
+S_n^{(\varepsilon)}(\alpha,z)}
+\frac{C\varkappa r(\alpha, z)}{ {\varphi(\sqrt{np_n})}
v^{3}|\alpha+S_n^{(\varepsilon)}(\alpha,z)|}.
\end{equation}
Analogously we get
%
%
\begin{equation}\label{ne9}
V(\alpha,z)=-\frac{\overline zS_n^{(\varepsilon)}(\alpha,z)}{\alpha
+S_n^{(\varepsilon)}(\alpha,z)}
+\frac{Cr(\alpha, z)}{\varphi(\sqrt{np_n}) v^{3}|\alpha
+S_n^{(\varepsilon)}(\alpha,z)|}.
\end{equation}
Inserting (\ref{ne8}) and (\ref{ne9}) in (\ref{new12}), we get
%
%
\begin{equation}\label{ne10}
\bigl(S_n^{(\varepsilon)}(\alpha,z)\bigr)^2+\alpha S_n^{(\varepsilon)}(\alpha
,z)+1-\frac{|z|^2S_n^{(\varepsilon)}(\alpha,z)}
{\alpha+S_n^{(\varepsilon)}(\alpha,z)}=\delta_n(z),
\end{equation}
where
\[
|\delta_n(\alpha,z)|\le\frac{C\varkappa}{\varphi(\sqrt{np_n})
v^{3} |S_n^{(\varepsilon)}(\alpha,z)+\alpha|}
\]
or equivalently
%
%
\begin{eqnarray}\label{99}
&&S_n^{(\varepsilon)}(\alpha,z) \bigl(\alpha+S_n^{(\varepsilon
)}(\alpha,z) \bigr)^2\nonumber\\[-8pt]\\[-8pt]
&&\qquad{}
+ \bigl(\alpha+S_n^{(\varepsilon)}(\alpha,z)
\bigr)-|z|^2S_n^{(\varepsilon)}(\alpha,z)=
\widetilde\delta_n(\alpha,z),\nonumber
\end{eqnarray}
where
$\widetilde\delta_n(\alpha,z)=\theta\frac{C\varkappa r(\alpha
,z)}{ \varphi(\sqrt{np_n}) v^{3} }$.

Furthermore, we introduce the notation
%
%
\begin{eqnarray}
Q_n^{(\varepsilon)}(\alpha,z)&:=&\bigl(\alpha+S_n^{(\varepsilon)}(\alpha
,z)\bigr)^2-|z|^2 \quad\mbox{and}\nonumber\\
Q(\alpha,z)&:=&\bigl(\alpha+S(\alpha,z)\bigr)^2-|z|^2,\\
P(\alpha,z)&:=&\alpha+S(\alpha,z) \quad\mbox{and}\quad
P^{(\varepsilon)}(\alpha,z):=\alpha+S_n^{(\varepsilon)}(\alpha,z).\nonumber
\end{eqnarray}
We may rewrite the last equation as
%
%
\begin{equation}\label{95}
S_n^{(\varepsilon)}(\alpha,z)=-\frac{P_n^{(\varepsilon)}(\alpha
,z)}{Q_n^{(\varepsilon)}(\alpha,z)}
+\widehat\delta_n(\alpha,z),
\end{equation}
where
%
%
\begin{equation}\label{delta}
\widehat\delta_n(\alpha,z)=\frac{\widetilde\delta_n(\alpha,z)}
{Q_n^{(\varepsilon)}(\alpha,z)}.
\end{equation}
%

Furthermore, we prove the following simple lemma.
\begin{lem}\label{lem03}Let $\alpha=u+iv$, $v>0$. Let $S(\alpha, z)$
satisfy the equation
%
%
\begin{equation}\label{eq:1}
S(\alpha,z)=-\frac{P(\alpha,z)}{Q(\alpha,z)},
\end{equation}
and $\operatorname{Im}\{S(\alpha,z)\}>0$.
Then the inequality
\[
1-|S(\alpha,z)|^2-\frac{|z|^2|S(\alpha,z)|^2}{|\alpha+S(\alpha,z)|^2}
\ge\frac v{v+1}
\]
holds.
\end{lem}
\begin{pf}
For $\alpha=u+iv$ with $v>0$, the Stieltjes transform $S(\alpha, z)$
satisfies the following equation:
%
%
\begin{equation}\label{eq:100}
S(\alpha,z)=-\frac{P(\alpha,z)}{Q(\alpha,z)}.
\end{equation}
Comparing the imaginary parts of both sides of this equation, we get
%
%
\begin{equation}\label{eq:2}
\operatorname{Im}\{P(\alpha,z)\}=\operatorname{Im}\{P(\alpha,z)\}
\frac{|P(\alpha,z)|^2+|z|^2}{|Q(\alpha,z)|^2} +v.
\end{equation}

Equations (\ref{eq:1}) and (\ref{eq:2}) together imply
%
%
\begin{equation}\label{eq:3}
\operatorname{Im}\{\alpha+S(\alpha,z)\} \biggl(1-\frac{|P(\alpha
,z)|^2+|z|^2}{|Q(\alpha,z)|^2} \biggr) =v.
\end{equation}
Since $v>0$ and $\operatorname{Im}\{\alpha+S(\alpha,z)\}>0$, it
follows that
\[
1-\frac{|P(\alpha,z)|^2+|z|^2}{|Q(\alpha,z)|^2}=1-|S(\alpha
,z)|^2-\frac{|z|^2|S(\alpha,z)|^2}{|\alpha+S(\alpha,z)|^2}>0.
\]
In particular we have
\[
|S(\alpha,z)|\le1.
\]
Equality (\ref{eq:3}) and the last remark together imply
\[
1-\frac{|P(\alpha,z)|^2+|z|^2}{|Q(\alpha,z)|^2}
=\frac v{\operatorname{Im}\{P(\alpha,z)\}}\ge\frac v{v+1}.
\]
The proof is complete.
\end{pf}

%
To compare the functions $S(\alpha,z)$ and $S_n(\alpha,z)$ we prove:
\begin{lem}\label{lem030}
Let
\[
|\widehat\delta_n(\alpha,z)|\le\frac v2.
\]
Then the following inequality holds
\[
1-\frac{|P_n^{(\varepsilon)}(\alpha,z)|^2+|z|^2}{|Q_n^{(\varepsilon
)}(\alpha,z)|^2}\ge\frac
v{4}.
\]
\end{lem}
%
\begin{pf}
By the assumption, we have
\[
\operatorname{Im}\{\widehat\delta_n(\alpha, z)+\alpha\}>\frac v2.
\]
Repeating the arguments of Lemma \ref{lem03} completes the proof.
\end{pf}

The next lemma provides a bound for the distance between the Stieltjes
transforms $S(\alpha,z)$ and $S_n^{(\varepsilon)}(\alpha, z)$.
\begin{lem}\label{lem04}
Let
\[
|\widehat\delta_n(\alpha,z)|\le\frac v8.
\]
Then
\[
\bigl|S_n^{(\varepsilon)}(\alpha,z)-S(\alpha,z)\bigr|\le\frac{4|\widehat
\delta_n(\alpha,z)|}v.
\]
\end{lem}
%
\begin{pf} Note that $S(\alpha,z)$ and $S_n^{(\varepsilon)}(\alpha
,z)$ satisfy the
equations
%
%
\begin{equation}\label{eq:10}
S(\alpha,z)=-\frac{P(\alpha,z)}{Q(\alpha,z)}
\end{equation}
and
%
%
\begin{equation}\label{eq:110}
S_n^{(\varepsilon)}(\alpha,z)=-\frac{P_n^{(\varepsilon)}(\alpha
,z)}{Q_n^{(\varepsilon)}(\alpha,z)}+\widehat\delta_n(\alpha,z),
\end{equation}
respectively.
These equations together imply
%
%
\begin{eqnarray}\quad
&&
S(\alpha,z)-S_n^{(\varepsilon)}(\alpha,z)\nonumber\\[-8pt]\\[-8pt]
&&\qquad=\frac{(S(\alpha
,z)-S_n^{(\varepsilon)}(\alpha,z))
(P_n^{(\varepsilon)}(\alpha,z)P(\alpha,z)+|z|^2)}
{Q(\alpha,z)Q_n^{(\varepsilon)}(\alpha,z)}+\widehat\delta_n(\alpha,z).\nonumber
\end{eqnarray}
Applying inequality $|ab|\le\frac12{(a^2+b^2)}$, we get
\begin{eqnarray*}
\biggl|1-\frac{P_n^{(\varepsilon)}(\alpha,z)P(\alpha,z)+|z|^2}
{Q(\alpha,z)Q_n^{(\varepsilon)}(\alpha,z)} \biggr|
&\ge&\frac12 \biggl(1-\frac{|P_n^{(\varepsilon)}(\alpha
,z)|^2+|z|^2}{|Q_n^{(\varepsilon)}(\alpha,z)|^2} \biggr)\\
&&{} +
\frac12 \biggl(1-\frac{|P(\alpha,z)|^2+|z|^2}{|Q(\alpha,z)|^2} \biggr).
\end{eqnarray*}
The last inequality and Lemmas \ref{lem03} and \ref{lem030}
together imply
\[
\biggl|1-\frac{P_n^{(\varepsilon)}(\alpha,z)P(\alpha,z)+|z|^2}
{Q(\alpha,z)Q_n^{\varepsilon)}(\alpha,z)} \biggr| \ge
\frac v4.
\]
This completes the proof of the lemma.
\end{pf}

%
To bound the distance between the distribution function $F_n(x,z)$ and the
distribution function $F(x,z)$
corresponding the Stieltjes transforms $S_n(\alpha,z)$ and $S(\alpha
,z)$ we use Corollary 2.3
from \cite{GT03}.
In the next lemma we give an integral bound for the distance between the
Stieltjes transforms $S(\alpha,z)$ and $S_n^{(\varepsilon)}(\alpha,z)$.
\begin{lem}
For $v\ge v_0(n)=c(\varphi(\sqrt{np_n}))^{-1/6}$ the inequality
\[
\int_{-\infty}^{\infty}\bigl|S(\alpha,z)-S_n^{(\varepsilon)}(\alpha
,z)\bigr|\,du\le
\frac{C(1+|z|^2)\varkappa}{\varphi(\sqrt{np_n}) v^7}
\]
holds.
\end{lem}
%
\begin{pf}Note that
%
%
\begin{equation}
\bigl|Q_n^{(\varepsilon)}\bigr|\ge\bigl|P_n^{(\varepsilon)}(\alpha
,z)-|z|\bigr|\bigl|P_n^{(\varepsilon)}(\alpha,z)+|z|\bigr|\ge v^2.
\end{equation}
It follows from here that $|\widehat\delta_n(\alpha, z)|\le\frac
{C}{v^5\varphi(\sqrt{np_n})}$ and
\[
|\widehat\delta_n(\alpha, z)|\le v/8
\]
for $v\ge c(\varphi(\sqrt{np_n}))^{-1/6}$.
Lemma \ref{lem04} implies that it is enough to prove the inequality
\[
\int_{-\infty}^{\infty}|\widehat\delta_n(\alpha, z)|\,du\le C\gamma_n,
\]
where $\gamma_n=\frac{C}{v^6\varphi(\sqrt{np_n})}$. By definition of
$\widehat\delta(\alpha,z)$, we have
%
%
\begin{equation}\label{np3}
\int_{-\infty}^{\infty}|\widehat\delta_n(\alpha, z)|\,du\le
\frac{c\varkappa}{v^3\varphi(\sqrt{np_n})}\int_{-\infty}^{\infty}
\frac{du}{|Q_n^{(\varepsilon)}(\alpha,z))|}.
\end{equation}
Furthermore, representation (\ref{95}) implies that
%
%
\begin{equation}
\frac{1}{|Q_n^{(\varepsilon)}(\alpha,z)|}\le\frac
{|S_n^{(\varepsilon)}(\alpha,z)|}{|P_n^{(\varepsilon)}(\alpha,z)|}
+\frac{|\widehat\delta_n(\alpha,z)|}{|P_n^{(\varepsilon)}(\alpha,z)|}.
\end{equation}
Note that, according to relation (\ref{ne10}),
%
%
\begin{equation}\label{012}
\frac1{|P_n^{(\varepsilon)}(\alpha,z)|}\le\frac
{|z|^2|S_n^{(\varepsilon)}(\alpha,z)|}
{|P_n^{(\varepsilon)}(\alpha,z)|^2}+\bigl|S_n^{(\varepsilon)}(\alpha,z)\bigr|+
\frac{|\delta_n(\alpha,z)|}{|P_n^{(\varepsilon)}(\alpha,z)|^2}.
\end{equation}
This inequality implies
%
%
\begin{eqnarray}\label{007}
\int_{-\infty}^{\infty}\frac{|S_n^{(\varepsilon)}(\alpha
,z)|}{|P_n^{(\varepsilon)}(\alpha,z)|}\,du&\le&
\frac{C(1+|z|^2)}{v^2}\int_{-\infty}^{\infty}\bigl|S_n^{(\varepsilon
)}(\alpha,z)\bigr|^2\,du\nonumber\\[-8pt]\\[-8pt]
&&{} +
\int_{-\infty}^{\infty}|\delta_n(\alpha,z)|
\frac{|S_n^{(\varepsilon)}(\alpha,z)|}{|P_n^{(\varepsilon)}(\alpha,z)|}\,du.\nonumber
\end{eqnarray}
It follows from relation (\ref{ne10}) that for $v>c(\varphi(\sqrt
{np_n}))^{-1/6}$,
%
%
\begin{equation}
|\delta_n(\alpha,z)|\le\frac{C\varkappa}{(\varphi(\sqrt{np_n}))v^4}<1/2.
\end{equation}
The last two inequalities together imply that for sufficiently large $n$
and $v>c(\varphi(\sqrt{np_n}))^{-1/6}$,
%
%
\begin{equation}\qquad\quad
\int_{-\infty}^{\infty}\frac{|S_n^{(\varepsilon)}(\alpha
,z)|}{|P_n^{(\varepsilon)}(\alpha,z)|}\,du\le
\frac
{C(1+|z|^2)}{v^2}\int_{-\infty}^{\infty}\bigl|S_n^{(\varepsilon)}(\alpha,z)\bigr|^2\,du
\le\frac{C(1+|z|^2)}{v^{3}}.
\end{equation}
Inequalities
(\ref{012}), (\ref{np3}) and the definition of
$\widehat\delta_n(\alpha,z)$ together imply
%
%
\begin{equation}\quad
\int_{-\infty}^{\infty}|\widehat\delta_n(\alpha,
z)|\,du\le\frac{C(1+|z|^2)}{ v^6\varphi(\sqrt{np_n})}
+\frac{C\varkappa}{v^4\varphi(\sqrt{np_n})}\int_{-\infty}^{\infty
}|\widehat\delta_n(\alpha,
z)|\,du.
\end{equation}
If we choose $v$ such that $\frac{C\varkappa}{v^4\varphi(\sqrt
{np_n}) }<\frac12$ we obtain
%
%
\begin{equation}
\int_{-\infty}^{\infty}|\widehat\delta_n(\alpha, z)|\,du\le
\frac{C(1+|z|^2)}{\varphi(\sqrt{np_n})v^6}.
\end{equation}
\upqed\end{pf}

In Section \ref{property} we show that the measure $\widetilde{\nu
}(\cdot,z)$ has
bounded support and bounded density for any $z$. To bound the distance
between the distribution functions $ \widetilde{F}_n^{(\varepsilon
)}(x,z)$ and $\widetilde F(x,z)$ we may apply
Corollary 3.2 from \cite{GT03} (see also Lemma~\ref{ap3} in the
\hyperref[app]{Appendix}). We take $V=1$ and
$v_0=C(\varphi(\sqrt{np_n}))^{-1/6}$. Then Lemmas \ref{lem03}
and \ref{lem030}
together imply
%
%
\begin{equation}\label{supremum}
\sup_x\bigl|F_n^{(\varepsilon)}(x,z)-F(x,z)\bigr|\le C\bigl(\varphi\bigl(\sqrt
{np_n}\bigr)\bigr)^{-1/6}.
\end{equation}
\upqed\end{pf}

\section{Properties of the measure $\widetilde{\nu}(\cdot,z)$}\label{property}

In this section we investigate the properties of the measure
$\widetilde{\nu}(\cdot,z)$. At
first note that there exists a solution $S(\alpha,z)$ of the equation
%
%
\begin{equation}\label{77}
S(\alpha,z)=-\frac{S(\alpha,z)+\alpha}{(S(\alpha,z)+\alpha)^2-|z|^2}
\end{equation}
such that, for $v>0$,
\[
\operatorname{Im}\{S(\alpha,z)\}\ge0
\]
and $S(\alpha,z)$ is an analytic function in the upper half-plane
$\alpha=u+iv$, $v>0$. This follows from the relative compactness of
the sequence of
analytic functions $S_n(\alpha,z)$, $n\in\mathbb N$. From
(\ref{eq:10}) it follows
immediately that
%
%
\begin{equation}
|S(\alpha,z)|\le1.
\end{equation}
Set $y=S(x,z)+x$ and consider equation (\ref{eq:10}) on the real line
%
%
\begin{equation}
y=-\frac{y}{y^2-|z|^2}+x
\end{equation}
or
%
%
\begin{equation}\label{eq:11}
y^3-xy^2+(1-|z|^2)y+x|z|^2=0.
\end{equation}

Set
%
%
\begin{eqnarray}
x_1^2&=&\frac{5+2|z|^2}{2}+\frac{(1+8|z|^2)^{3/2}-1}{8|z|^2
},\nonumber\\[-8pt]\\[-8pt]
x_2^2&=&\frac{5+2|z|^2}{2}-\frac{(1+8|z|^2)^{3/2}+1}{8|z|^2 }.\nonumber
\end{eqnarray}
It is straightforward to check that
$\sqrt{3(1-|z|^2)}\le|x_1|$
and
$x_2^2<0$ for $|z|<1$ and
$x_2^2=0$ for $|z|=1$, and $x_2^2>0$ for $|z|>1$.
\begin{lem}\label{lem01}
In the case $|z|\le1$ equation (\ref{eq:11}) has one real root for
$|x|\le|x_1|$
and three real roots for
$|x|>|x_1|$.
In the case $|z|>1$ equation (\ref{eq:11}) has one real root for
$|x_2|\le x\le|x_1|$
and three real roots for
$|x|\le|x_2|$
or for
$|x|\ge|x_1|$.
\end{lem}
\begin{pf}
Set
\[
L(y):=y^3-xy^2+(1-|z|^2)y+x|z|^2.
\]
We consider the roots of the equation
%
%
\begin{equation}
L'(y)=3y^2-2xy+(1-|z|^2)=0.
\end{equation}
The roots of this equation are
\[
y_{1,2}=\frac{x\pm\sqrt{x^2-3(1-|z|^2)}}3.
\]
This implies that, for $|z|\le1$ and for
\[
|x|\le\sqrt{3(1-|z|^2)}
\]
equation (\ref{eq:11}) has one real root.
Furthermore, direct calculations show that
\[
L(y_1)L(y_2)=\tfrac1{27} \bigl(-4|z|^2x^4+(8|z|^4+20|z|^2-1)x^2
+4(1-|z|^2)^3 \bigr).
\]

Solving the equation $L(y_1)L(y_2)=0$ with respect to $x$, we get
for $|z|\le1$ and $\sqrt{3(1-|z|^2)}\le|x|\le|x_1|$
\[
L(y_1)L(y_2)\ge0,
\]
and
for $|z|\le1$ and $|x|>\sqrt{\frac{20+8|z|^2}{8}
+\frac{(1+8|z|^2)^{3/2}-1}
{8|z|^2 }}$
\[
L(y_1)L(y_2)< 0.
\]
These relations imply that
for $|z|\le1$ the function
$L(y)$ has three real roots for $|x|\ge|x_1|$
and one real root for $|x|<|x_1|$.

Consider the case $|z|>1$ now. In this case $y_{1,2}$ are real for
all $x$ and $x_2^2>0$.
Note that
\[
L(y_1)L(y_2)\le0
\]
for $|x|\le|x_2|$ and for $|x|\ge|x_1|$ and
\[
L(y_1)L(y_2)>0
\]
for $|x_2|<x<|x_1|$. These implies that for $|z|>1$ and for
$|x_2|<x<|x_1|$ the function $L(y)$ has one real root and for
$|x|\le|x_2|$ or for $|x|\ge|x_1|$ the function $L(y)$ has three real
roots. The lemma is proved.
\end{pf}
\begin{rem}\label{rem:01}
From Lemma \ref{lem01} it follows that the
measure $\widetilde{\nu}(x,z)$ has
a density $p(x,z)=\lim_{v\to0}\operatorname{Im}{S(\alpha,z)}$ and:
\begin{itemize}
\item{$p(x,z)\le1$, for all $x$ and $z$;}
\item{for $|z|\le1$, if $|x|\ge x_1$, then $p(x,z)=0$;}
\item{for $|z|\ge1$, if $|x|\ge x_1$ or $|x|\le x_2$, then $p(x,z)=0$;}
\item{$p(x,z)>0$ otherwise.}
\end{itemize}
\end{rem}

Introduce the function
%
%
\begin{equation}
g(s,t):=\cases{\dfrac{2s}{s^2+t^2}, &\quad if $s^2+t^2>1$,\vspace*{2pt}\cr
2s, &\quad otherwise.}
\end{equation}
It is well known that for $z=s+it$ the logarithmic potential of uniform
distribution on the unit disc is
%
%
\begin{equation}\label{potential}\qquad
U_0(z):=\iint\ln{\frac1{|z-x+iy|}}\,dG(x,y)=\cases{\dfrac
12(1-|z|^2), &\quad if $|z|\le1$,\cr
-\ln|z|, &\quad if $|z|>1$,}
\end{equation}
and
%
%
\begin{equation}
\frac{ \partial}{ \partial s}\iint\ln{\frac
1{|z-x+iy|}}\,dG(x,y)=-\frac12g(s,t).
\end{equation}

According to Lemma 4.4 in Bai \cite{Bai1997}, we have, for $z=s+it$,
%
%
\begin{equation}\label{partial}
\frac{ \partial}{ \partial s} \biggl(\int_0^{\infty}\log x\nu
(dx,z) \biggr)
=\frac12g(s,t).
\end{equation}
According to Remark \ref{rem:01}, we have, for $|z|\ge1$,
%
%
\begin{equation}
\ln(|x_2|/|z|)\le U_{\widetilde{\nu}}(z)+\ln|z|\le\ln(|x_1|/|z|).
\end{equation}
This implies that
%
%
\begin{equation}\label{qq}
{\lim_{|z|\to\infty}}|U_{\widetilde{\nu}}(z)- U_0(z)|=0.
\end{equation}
Since
%
%
\begin{equation}
\int_{-\infty}^{\infty}\log|x| \widetilde{\nu}(dx,z)=\int
_{0}^{\infty}\log x {\nu}(dx,z)
\end{equation}
we get
%
%
\begin{equation}
\frac{ \partial}{ \partial s} \biggl(\int_{-\infty}^{\infty}\log
|x| \widetilde{\nu}(dx,z) \biggr)
=\frac12g(s,t).
\end{equation}
Comparing equalities (\ref{partial}) and (\ref{potential}) and using
relation (\ref{qq}), we obtain
%
%
\begin{equation}\label{main}
U_0(z)=-\int_{0}^{\infty}\ln x \nu(dx, z)=-\int_{-\infty
}^{\infty}\ln|x| \widetilde{\nu}(dx,z)=U_{\mu}(z).
\end{equation}

\section{The smallest singular value}\label{singular}

Let $\mathbf X^{(\varepsilon)}=\frac1{\sqrt{np_n}} (\varepsilon
_{jk}X_{jk} )_{j,k=1}^n$ be an $n\times n$ matrix
with independent entries $\varepsilon_{jk}X_{jk}$, $j,k=1,\ldots,n$. Assume
that $\mathbf{E}X_{jk}=0$ and $\mathbf{E}X_{jk}^2=1$ and let
$\varepsilon_{jk}$ denote Bernoulli random variables with $p_n=
\Pr\{\varepsilon_{jk}=1\}$, $j,k=1,\ldots,n$.
Denote by $s_1^{(\varepsilon)}(z)\ge\cdots\ge s_n^{(\varepsilon
)}(z)$ the
singular values of the matrix $\mathbf X^{(\varepsilon)}(z):=\mathbf
X^{(\varepsilon)}-z\mathbf I$. In this section we prove a bound for
the minimal singular
value of the matrices $\mathbf X^{(\varepsilon)}(z)$.
We prove the following result.
%
\begin{theorem}\label{thm1}Let $X_{jk}, j,k \in\mathbf{N}$,\vspace*{2pt} be
independent random complex variables with $\mathbf{E}X_{jk}=0$ and
$\mathbf{E}|X_{jk}|^2=1$, which are
uniformly integrable, that is,
%
%
\begin{equation}\label{unif0}
\sup_{j,k}\mathbf{E}|X_{jk}|^2I_{\{|X_{jk}|>M\}}\to0 \qquad\mbox
{as } M\to\infty.
\end{equation}
Let $\varepsilon_{jk}$, $j,k=1,\ldots,n$, be independent Bernoulli
random variables with
$p_n:=\Pr\{\varepsilon_{jk}=1\}$.
Assume that $\varepsilon_{jk}$ are independent from $X_{jk}, j,k\in
\mathbf{N}$, in aggregate. Let $p_n^{-1}=\mathcal O(n^{1-\theta})$ for
some $0<\theta\le1$.
Let $K\ge1$. Then there exist constants $c, C, B >0$ depending on
$\theta$ and $K$ such that for any $z\in\mathbb C$ and positive
$\varepsilon$ we have
%
%
\begin{equation}
\Pr\bigl\{s_n^{(\varepsilon)}(z)\le\varepsilon/n^B ; s_1^{(\varepsilon
)}(z)\le Kn\sqrt{p_n}\bigr\}\le
\exp\{ - c p_n n \}+\frac{C\sqrt{\ln n}}{\sqrt{np_n}}.
\end{equation}
\end{theorem}
%
\begin{rem}\label{thm11} Let $X_{jk}$ be i.i.d. random variables
with $\mathbf{E}X_{jk}=0$ and $\mathbf{E}|X_{jk}|^2=1$. Then
condition (\ref{unif0}) holds.
\end{rem}
\begin{rem}
Consider the event $A$ that there exists at least one row
with zero entries only. Its probability is given by
%
%
\begin{equation}
\Pr\{A\}\ge1-\bigl(1-(1-p_n)^n\bigr)^n.
\end{equation}
Simple calculations show that if $np_n\le\ln n$ for all $n\ge1$, then
%
%
\begin{equation}
\Pr\{A\}\ge\delta>0.
\end{equation}
Hence in the case $np_n\le\ln n$ and $np_n\to\infty$ we have no
invertibility with positive probability.
\end{rem}
%
\begin{rem}
The proof of Theorem \ref{thm1} uses ideas of Rudelson and Vershynin
\cite{RV},
to classify with high probability vectors $\mathbf x$ in the
$(n-1)$-dimensional unit
sphere $\mathcal{S}^{n-1}$ such that $\Vert\mathbf X^{(\varepsilon)}(z)
\mathbf x\Vert_2$ is
extremely small into two classes, called compressible and
incompressible vectors.

We develop our approach for shifted sparse and normalized matrices
$\mathbf X^{(\varepsilon)}(z)$.
The generalization to the case of complex sparse and shifted matrices
$\mathbf X^{(\varepsilon)}(z)$ is straightforward.
For details see, for example, the paper of G\"{o}tze and Tikhomirov
\cite{GT07} and the proof of the Lemma \ref{lem02} below.
\end{rem}
%
%
\begin{rem}
We may relax the condition $p_n^{-1}=\mathcal O(n^{1-\theta})$ to
$p_n^{-1}= o(n/\break\ln^2{n})$. The quantity $B$ in Theorem \ref{thm1}
should be of order $\ln n$ in this case. See Remark \ref{rem:relax}
for details.
\end{rem}
\begin{lem}\label{lem02}Let $\mathbf x=(x_1,\ldots,x_n)\in\mathcal
S^{n-1}$ be a fixed unit vector and $\mathbf X^{(\varepsilon)}(z)$ be
a matrix as in Theorem \ref{thm1}.
Then there exist some positive absolute constants $\gamma_0$ and $c_0$ such
that for any $0<\tau\le\gamma_0 $
%
%
\begin{equation}\label{imp1}
\Pr\bigl\{\bigl\|\mathbf X^{(\varepsilon)}(z) \mathbf x\bigr\|_2\le\tau\bigr\}\le\exp
\{-c_0 n p_n\}.
\end{equation}
\end{lem}
%
\begin{pf}
Recall that $\mathbf{E}X_{ij}=0$ and $\mathbf{E}|X_{ij}|^2=1$. Assume
first that $X_{ij}$ are real independent r.v. with
mean zero, and variance at least $1$. Let $X^{(\varepsilon)}_{ij}=
X_{ij} \varepsilon_{ij}$ with
independent Bernoulli variables which are independent of $X_{ij}$ in
aggregate and let $z=0$.
Assume also that $\mathbf x$ is a real vector. Then
%
%
\begin{equation}
\bigl\Vert\mathbf X^{(\varepsilon)} \mathbf x\bigr\Vert_2^2 = \frac1{np_n}\sum
_{j=1}^n \Biggl|\sum_{k=1}^nx_kX_{jk}\varepsilon_{jk} \Biggr|^2=:\frac
1{np_n}\sum_{k=1}^n\zeta_j^2.
\end{equation}
By Chebyshev's inequality we have
%
%
\begin{eqnarray}\label{cheb1}
\Pr\Biggl\{\sum_{j=1}^n\zeta_j^2<\tau^2np_n\Biggr\}&=&\Pr\Biggl\{\frac{\tau
^2np_n}2-\frac12\sum_{j=1}^n\zeta_j^2>0\Biggr\}\nonumber\\[-8pt]\\[-8pt]
&\le&
\exp\{np_n\tau^2t^2/2\}\prod_{j=1}^n\mathbf{E}\exp\{-{t^2\zeta
_j^2/2}\}.\nonumber
\end{eqnarray}
Using $e^{-t^2/2}=\mathbf{E}\exp\{it\xi\}$, where $\xi$ is a
standard Gaussian random variable, we obtain
%
%
\begin{eqnarray}\label{4.7}
&&\Pr\Biggl\{\sum_{j=1}^n\zeta_j^2<\tau^2np_n\Biggr\}\nonumber\\[-8pt]\\[-8pt]
&&\qquad\le\exp\{np_n\tau^2t^2/2\}
\prod_{j=1}^n\mathbf{E}_{\xi_j}\prod_{k=1}^n
\mathbf{E}_{\varepsilon_{jk}X_{jk}}\exp\{it\xi_jx_k\varepsilon
_{jk}X_{jk}\},\nonumber
\end{eqnarray}
where $\xi_j$, $j=1,\ldots,n$, denote i.i.d. standard Gaussian
r.v.s and $\mathbf{E}_Z$ denotes expectation with respect to $Z$
conditional on
all other r.v.s.
For every $\alpha, x\in[0,1]$ and $\rho\in(0,1)$ the following
inequality holds:
%
%
\begin{equation}\label{goetze}
\alpha x+1-\alpha\le x^\beta\vee\biggl(\frac{\rho}{\alpha}
\biggr)^{{\beta}/({1-\beta})}
\end{equation}
(see \cite{Bickel},
inequality (3.7)). Take $\alpha=\Pr\{|\xi_j|\le C_1\}$
for some absolute positive constant $C_1$ which will be chosen later.
Then it follows from (\ref{4.7}) that
%
%
\begin{eqnarray}\label{i0}
&&\Pr\Biggl\{\sum_{j=1}^n\zeta_j^2<\tau^2np_n\Biggr\}\nonumber\\
&&\qquad\le\exp\{np_n\tau
^2t^2/2\}\\
&&\qquad\quad{} \times\prod_{j=1}^n
\Biggl(\alpha\Biggl|\mathbf{E}_{\xi_j} \Biggl(\prod_{k=1}^n\mathbf
{E}_{\varepsilon_{jk}X_{jk}}\exp\{it\xi_jx_k\varepsilon_{jk}X_{jk}\}
\Big||\xi_j|\le C_1 \Biggr) \Biggr|+1-\alpha\Biggr).\hspace*{-19pt}\nonumber
\end{eqnarray}
Furthermore, we note that
%
%
\begin{eqnarray}\label{4.8}
&&|\mathbf{E}_{\varepsilon_{jk}X_{jk}}\exp\{it\xi_jx_k\varepsilon
_{jk}X_{jk}\}|\nonumber\\
&&\qquad\le\exp\biggl\{\frac12(|\mathbf{E}_{\varepsilon_{jk}X_{jk}}\exp\{it\xi
_jx_k\varepsilon_{jk}X_{jk}\}|^2-1)\biggr\}\nonumber\\[-8pt]\\[-8pt]
&&\qquad
\le\exp\biggl\{-{p_n} \biggl((1-p_n)\bigl(1-\operatorname{Re}{f_{jk}(tx_k\xi
_j)}\bigr)\nonumber\\
&&\qquad\quad\hspace*{62.8pt}{}+\frac{p_n}2\bigl(1-|f_{jk}(tx_k\xi_j)|^2\bigr) \biggr)\biggr\},\nonumber
\end{eqnarray}
where $f_{jk}(u)=\mathbf{E}\exp\{iuX_{jk}\}$.
Assuming (\ref{unif0}), choose a constant $M>0$ such that
%
%
\begin{equation}
\sup_{jk}\mathbf{E}|X_{jk}|^2I_{\{|X_{jk}|>M\}}\le1/2.
\end{equation}
Since $1-\cos x\ge11/24 x^2$ for $|x|\le1$, conditioning on the event
$|\xi_j|\le C_1$, we get for
$0<t \le1/ (MC_1)$
%
%
\begin{eqnarray}\label{i11}
1-\operatorname{Re}f_{jk}(tx_k\xi_j)&=&\mathbf{E}_{X_{jk}}\bigl(1-\cos
(tx_kX_{jk}\xi_j)\bigr)\nonumber\\[-8pt]\\[-8pt]
&\ge&\tfrac{11}{24}{t^2x_k^2\xi_j^2}\mathbf
{E}|X_{jk}|^2I_{\{|X_{jk}|\le M\}},\nonumber
\end{eqnarray}
and similarly
%
%
\begin{eqnarray}\label{i1}
1-|f_{jk}(tx_k\xi_j)|^2 &=&\mathbf{E}_{X_{jk}}\bigl(1-\cos(tx_k\widetilde
X_{jk}\xi_j)\bigr)\nonumber\\[-8pt]\\[-8pt]
&\ge&\tfrac{11}{24}{t^2x_k^2\xi_j^2}
\mathbf{E}|\widetilde X_{jk}|^2I_{\{|X_{jk}|\le
M\}}.\nonumber
\end{eqnarray}

It follows from (\ref{4.8}) for $0<t<1/{(MC_1)}$ and for some constant $c>0$
%
%
\begin{equation}\label{i2}
|\mathbf{E}_{\varepsilon_{jk}X_{jk}}\exp\{it\xi_jx_k\varepsilon
_{jk}X_{jk}\}|\le\exp\{-{cp_n}t^2x_k^2\xi_j^2\}.
\end{equation}
This implies that conditionally on $|\xi_j|\le C_1$ and for $ 0<t\le1/(MC_1)$
%
%
\begin{equation}\label{i3}
\Biggl|\prod_{k=1}^n\mathbf{E}_{\varepsilon_{jk}X_{jk}}\exp\{it\xi
_jx_k\varepsilon_{jk}X_{jk}\}\Biggr|\le\exp\{-{cp_nt^2\xi_j^2}\}.
\end{equation}
Let $\Phi_0(x):=2\Phi(x)-1$, $x>0$, where $\Phi(x)$ denotes the
standard Gaussian distribution function.
It is straightforward to show that
%
%
\begin{eqnarray}
&&\mathbf{E}_{\xi_j} (\exp\{-cp_nt^2\xi_j^2\} ||\xi_j|\le
C_1 )
\nonumber\\[-8pt]\\[-8pt]
&&\qquad=\frac1{\sqrt{1+2ct^2p_n}}
\frac{\Phi_0 (C_1\sqrt{1+2t^2cp_n} )}{\Phi_0(C_1)}.\nonumber
\end{eqnarray}
We may choose $C_1$ large enough such that following inequalities hold:
%
%
\begin{equation}\label{4.9}
\mathbf{E}_{\xi_j} (\exp\{-cp_nt^2\xi_j^2\} ||\xi_j|\le
C_1 )\le
\exp\{-ct^2p_n/24\}
\end{equation}
for all $|t|\le1/(MC_1)$.
Inequalities (\ref{4.7}), (\ref{goetze}), (\ref{4.8}), (\ref{4.9})
together imply that for any $\beta\in(0,1)$
%
%
\begin{eqnarray}
&&\Pr\Biggl\{\sum_{j=1}^n\zeta_j^2<\tau^2np_n \Biggr\}\nonumber\\[-8pt]\\[-8pt]
&&\qquad\le\exp\{
np_n\tau^2t^2/2\}
\biggl(\exp\{-c\beta n t^2p_n/24\}+ \biggl(\frac{\beta}{\alpha}
\biggr)^{{n\beta}/({1-\beta})} \biggr).\nonumber
\end{eqnarray}
Without loss of generality we may take $C_1$ sufficiently large, such
that $\alpha\ge4/5$ and choose
$\beta=2/5$. Then we obtain
%
%
\begin{eqnarray}
&&\Pr\Biggl\{\sum_{j=1}^n\zeta_j^2<\tau^2np_n\Biggr\}
\nonumber\\[-8pt]\\[-8pt]
&&\qquad\le\exp\{np_n\tau^2t^2/2\}
\biggl(\exp\{-ct^2np_n/60\}
+ \biggl(\frac12 \biggr)^{{2n}/3} \biggr).\nonumber
\end{eqnarray}
For $\tau<\frac{\sqrt c}{\sqrt{60}}$ we conclude from here that for
$|t|\le1/(MC_1)$
%
%
\begin{equation}\label{i10}
\Pr\Biggl\{\sum_{j=1}^n\zeta_j^2<\tau^2np_n\Biggr\}\le\exp\{-{ct^2np_n/120}\}.
\end{equation}
Inequality (\ref{i10}) implies that inequality (\ref{imp1}) holds
with some positive constant $c_0>0$.
This completes the proof in the real case.

Consider now the general case.
Let $X_{jk}=\xi_{jk}+i\eta_{jk}$ with $i=\sqrt{-1}$ with $\mathbf
{E}|X_{jk}|^2=1$ and $x_k=u_k+iv_k$ and $z=u+iv$.
In this notation we have
%
%
\begin{eqnarray}\label{r1}
&&\Pr\bigl\{\bigl\|\bigl(\mathbf X^{(\varepsilon)}-z\mathbf I\bigr)\mathbf x\bigr\|_2\le\tau\bigr\}
\nonumber\\
&&\qquad\le
\exp\{\tau^2 np_nt^2/2\}\nonumber\\
&&\qquad\quad{}\times
\min\Biggl\{\mathbf{E}\exp\Biggl\{-t^2\sum_{j=1}^n \Biggl|\sum
_{k=1}^n(\xi_{jk}u_k-\eta_{jk}v_k)\varepsilon_{jk}\nonumber\\[-8pt]\\[-8pt]
&&\hspace*{145pt}{}-\sqrt{np_n}(uu_j-vv_j) \Biggr|^2\Big/2
\Biggr\},\nonumber\\
&&\qquad\quad\hspace*{37.1pt} \mathbf{E}\exp\Biggl\{-t^2\sum_{j=1}^n
\Biggl|\sum_{k=1}^n(\xi_{jk}v_k+\eta_{jk}u_k)\varepsilon_{jk}\nonumber\\
&&\hspace*{145pt}{}-\sqrt
{np_n}(vu_j+uv_j)
\Biggr|^2\Big/2 \Biggr\} \Biggr\}.\nonumber
\end{eqnarray}
Note that for $\mathbf x=(x_1,\ldots, x_n)\in S^{(n-1)}$ (the unit
sphere in $\mathbb C^n$) and for any set
$A\subset\{1,\ldots, n\}$
%
%
\begin{equation}\label{i4}
\max\biggl\{\sum_{k\in A}|x_k|^2, \sum_{k\in A^c}|x_k|^2\biggr\}\ge1/2.
\end{equation}
For any $j=1,\ldots, n$ we introduce the set $A_j$ as follows:
%
%
\begin{equation}
A_j:=\bigl\{k\in\{1,\ldots, n\}\dvtx\mathbf{E}|\xi_{jk}u_k-\eta
_{jk}v_{k}|^2\ge|x_k|^2/2\bigr\}.
\end{equation}
It is straightforward to check that for any $k\notin A_j$
%
%
\begin{equation}
\mathbf{E}|\eta_{jk}u_k+\xi_{jk}v_k|^2\ge|x_k|^2/2.
\end{equation}
According to inequality (\ref{i4}), for any $j=1,\ldots,n$, there
exists a set $B_j$ such that
%
%
\begin{equation}
\sum_{k\in B_j}|x_k|^2\ge1/2
\end{equation}
and for any $k\in B_j$
%
%
\begin{equation}\label{t1}
\mathbf{E}|\xi_{jk}u_k-\eta_{jk}v_{k}|^2\ge|x_k|^2/2
\end{equation}
or
%
%
\begin{equation}\label{t2}
\mathbf{E}|\eta_{jk}u_k+\xi_{jk}v_k|^2\ge|x_k|^2/2.
\end{equation}
Introduce the following random variables for any $j,k=1,\ldots,n$
%
%
\begin{equation}
\widetilde\zeta_{jk}:=\xi_{jk}u_k-\eta_{jk}v_{k}
\end{equation}
and
%
%
\begin{equation}
\widehat\zeta_{jk}:=\eta_{jk}u_k+\xi_{jk}v_k.
\end{equation}
Inequalities (\ref{t1}) and (\ref{t2}) together imply that one of the
following two inequalities
%
%
\begin{equation}\label{t11}
\operatorname{card} \{j\dvtx\mbox{for any } k\in B_j\mbox{ } \mathbf
{E}|\widehat\zeta_{jk}|^2\ge|x_k|^2/2 \}\ge n/2
\end{equation}
or
%
%
\begin{equation}\label{t21}
\operatorname{card} \{j\dvtx\mbox{for any } k\in B_j\mbox{ } \mathbf
{E}|\widetilde\zeta_{jk}|^2\ge|x_k|^2/2 \}\ge n/2
\end{equation}
holds.
If (\ref{t11}) holds we shall bound the first term on the right-hand
side of (\ref{r1}). In the other case we shall bound the second term.
In what follows we may repeat the arguments leading to inequalities
(\ref{i0})--(\ref{i3}).
Thus the lemma is proved.
\end{pf}

For any $q_n\in(0,1)$ and $K>0$ to be chosen later we define
$K_n:=Kn\sqrt{p_n}$, $\widehat q_n:=q_n/ (\ln(2/p_n)\ln K_n
)$ and $\widehat p_n:=p_n/ (\ln(2/p_n)\ln K_n )$. Without loss
of generality we shall assume that
%
%
\begin{equation}\label{kk}
\ln K_n/|\ln\gamma_0|\ge1 \quad\mbox{and}\quad \ln K_n>1.
\end{equation}
\begin{prop}\label{prop}
Assume there exist an absolute constant $c>0$ and values $\gamma_n,
q_n\in(0,1)$ such that for any $\mathbf x\in\mathcal C\subset
\mathcal S^{(n-1)}$
%
%
\begin{equation}\label{g02}
\Pr\bigl\{\bigl\|\mathbf X^{(\varepsilon)}(z)\mathbf x\bigr\|_2\le\gamma_n
\mbox{ and } \bigl\|\mathbf X^{(\varepsilon)}(z)\bigr\|\le K_n\bigr\}\le\exp\{
-c nq_n\}
\end{equation}
holds.
Then there exists a constant $\delta_0>0$ depending on $K$ and $c$
only such that, for $k<\delta_0 n \widehat q_n$,
\[
\Pr\Bigl\{\inf_{\mathbf x\in\mathcal S^{k-1}\cap\mathcal C}\bigl\|\mathbf
X^{(\varepsilon)}(z)\mathbf x\bigr\|_2\le\gamma_n/2\mbox{ and }\bigl\|\mathbf
X^{(\varepsilon)}(z)\bigr\|\le K_n\Bigr\}\le\exp\{-c nq_n/8\}.
\]
\end{prop}
\begin{pf}
Let $\eta>0$ to be chosen later. There exists an $\eta$-net
$\mathcal N$ in $\mathcal S^{k-1}\cap\mathcal C$ of cardinality
$|\mathcal N|\le(\frac{3}{\eta})^{2k}$
(see, e.g., Lemma 3.4 in \cite{rud06}).
By condition (\ref{g02}), we have for $\tau\le\gamma_n$
%
%
\begin{eqnarray}
&&\Pr\bigl\{\mbox{there exists }\mathbf x\in\mathcal N\dvtx\bigl\|\mathbf
X^{(\varepsilon)}(z)\mathbf x\bigr\|_2 < \tau\mbox{ and }\bigl\|\mathbf
X^{(\varepsilon)}(z)\bigr\|\le K_n\bigr\}\nonumber\\[-8pt]\\[-8pt]
&&\qquad\le
\biggl(\frac{3}{\eta}\biggr)^{2k}\exp\{-c n q_n\}.\nonumber
\end{eqnarray}
Let $V$ be the event that $\|\mathbf X^{(\varepsilon)}(z)\|\le K_n$
and $\|\mathbf X^{(\varepsilon)}(z) \mathbf y\|_2\le\frac12 \tau$
for some point \mbox{$\mathbf y\in\mathcal S^{(k-1)}\cap\mathcal C$}.
Assume that $V$ occurs and choose a point $\mathbf x\in\mathcal N$
such that $\|\mathbf y-\mathbf x\|_2\le\eta$.
Then
%
%
\begin{equation}\qquad
\bigl\|\mathbf X^{(\varepsilon)}(z)\mathbf x\bigr\|_2\le\bigl\|\mathbf
X^{(\varepsilon)}(z)\mathbf y\bigr\|_2+\bigl\|\mathbf X^{(\varepsilon)}(z)\bigr\|\|
\mathbf x-\mathbf y\|_2\le\tfrac12\tau+K_n\eta=\tau,
\end{equation}
if we set $\eta=\tau/(2K_n)$. Hence,
%
%
\begin{equation}
\Pr(V)\le\biggl( \biggl(\frac3{\eta} \biggr)^{2\delta_0/(\ln K_n\ln
(2/p_n))}\exp\biggl\{-\frac{c_0}4\biggr\} \biggr)^{nq_n}.
\end{equation}
Note that under assumption (\ref{kk}) we have
%
%
\begin{equation}
\frac{2\ln(3/\eta)}{\ln2\ln K_n}\le10.
\end{equation}
Choosing $\delta_0=\frac{c}{80}$ and $\tau=\gamma_n$, we complete
the proof.
\end{pf}
%


Following Rudelson and Vershynin \cite{RV}, we shall partition the
unit sphere $\mathcal S^{(n-1)}$ into the two sets of so-called
compressible and incompressible vectors, and we will show the
invertibility of $\mathbf X$ on each set separately.
\begin{defn}
Let $\delta,\rho\in(0,1)$. A vector $\mathbf x\in\mathbb R^n$ is
called \textit{sparse} if $|{\operatorname{supp}}(\mathbf x)|\le\delta n$.
A vector $\mathbf x\in\mathcal S^{(n-1)}$ is called \textit{compressible}
if $\mathbf x$ is within Euclidean distance $\rho$ from the set of all
sparse vectors.
A vector $\mathbf x\in\mathcal S^{(n-1)}$ is called
\textit{incompressible} if it is not compressible.
\end{defn}

The sets of sparse, compressible and incompressible vectors depending
on $\delta$ and $\rho$ will be denoted by
%
%
\begin{equation}
\mathit{Sparse}(\delta),\qquad \mathit{Comp}(\delta,
\rho),\qquad
\mathit{Incomp}(\delta, \rho),
\end{equation}
respectively.

%
\begin{lem}\label{lem031} Let $\mathbf X^{(\varepsilon)}(z)$ be a
random matrix as in Theorem \ref{sparse}, and let
$K_n=Kn\sqrt{p_n}$ with a constant $K\ge1$. Assume there exist an
absolute constant $c>0$ and values $\gamma_n, q_n\in(0,1)$ such that
for any $\mathbf x\in\mathcal C\subset\mathcal S^{(n-1)}$
%
%
\begin{equation}\label{g01}
\Pr\bigl\{\bigl\|\mathbf X^{(\varepsilon)}(z)\mathbf x\bigr\|_2\le\gamma_n
\mbox{ and }\bigl\|\mathbf X^{(\varepsilon)}(z)\bigr\|\le K_n\bigr\}\le\exp\{
-c nq_n\}
\end{equation}
holds.
Then there exist $\delta_1, c_1$ that depend
on $K$ and $c$ only, such that
%
%
\begin{eqnarray}
&&\Pr\Bigl\{\inf_{\mathbf x\in \mathrm{Comp}(\delta_1\widehat
q_n, \rho_{n})\cap\mathcal C}\bigl\|\mathbf X^{(\varepsilon)}(z) \mathbf
x\bigr\|_2\le\gamma_n \mbox{ and } \bigl\|\mathbf X^{(\varepsilon)}(z) \bigr\|\le
K_n \Bigr\}\nonumber\\[-8pt]\\[-8pt]
&&\qquad\le\exp\{-c_1 n q_n\},\nonumber
\end{eqnarray}
where $\rho_n:=\gamma_n/(4K_n)$.
\end{lem}
%
\begin{pf}At first we estimate the invertibility for sparse vectors.
Let $k=[\delta_1 n\widehat q_n]$ with some positive constant $\delta_1$
which will be chosen later.
According to Proposition \ref{prop} for any $\delta_1\le\delta_0$
and for any $\tau\le\gamma_n/2$, we have the following inequality:
\begin{eqnarray*}
&&\Pr \Bigl\{\inf_{\mathbf x\in\mathit{Sparse}(\delta_1\widehat
p_n)\cap\mathcal C}\bigl\|\mathbf X^{(\varepsilon)} (z)\mathbf x\bigr\|_2\le
\tau\mbox{ and } \bigl\|\mathbf X^{(\varepsilon)} (z)\bigr\|
\le K_n \Bigr\} \\
&&\qquad=
\Pr\Bigl\{\mbox{there exists } \sigma, |\sigma|=k\dvtx
\inf_{\mathbf x\in\mathbb R^{\sigma}\cap\mathcal C, \|\mathbf x\|
_2=1}\bigl\|\mathbf X ^{(\varepsilon)} (z)\mathbf x\bigr\|_2
\le\tau\\
&&\hspace*{190pt}\mbox{ and } \|\mathbf
X^{(\varepsilon)} (z)\|\le K_n \Bigr\}\\
&&\qquad\le\pmatrix{n\cr k}\exp\{-c_0nq_n/8\}.
\end{eqnarray*}
Using Stirling's formula, we get for some absolute positive constant $C$
%
%
\begin{equation}
\frac1n\ln\pmatrix{n \cr k}\le-C\delta_1 \widehat q_n\ln(\delta
\widehat q_n).
\end{equation}
We may choose $\delta_1$ small enough that
%
%
\begin{equation}
\frac1n\ln\pmatrix{n \cr k}\le c_0q_n/{16}.
\end{equation}
Thus we get
%
%
\begin{equation}\qquad
\Pr\Bigl\{\inf_{\mathbf x\in\mathit{Sparse}(\delta_1\widehat p_n)\cap
\mathcal C}\bigl\|\mathbf X ^{(\varepsilon)} (z)\mathbf
x\bigr\|_2\le\tau\mbox
{ and } \bigl\|\mathbf X^{(\varepsilon)} (z)\bigr\|
\le K_n \Bigr\}\le\exp\{-c_1 nq_n\}.\hspace*{-38pt}
\end{equation}
Choose $\rho:=\gamma:=\gamma_n/4$. Let $V$ be the event that $\|
\mathbf X^{(\varepsilon)}(z)\|\le K_n$ and $\|\mathbf X^{(\varepsilon
)}(z) \mathbf y\|_2\le\gamma_1$ for some point $\mathbf y\in
\mathit{Comp}(\delta_1\widehat p_n,\rho K_n^{-1})$.
Assume that $V$ occurs and choose a point $\mathbf x\in
\mathit{Sparse}(\delta_1\widehat p_n)$
such that $\|\mathbf y-\mathbf x\|_2\le\rho K_n^{-1}$.
Then
%
%
\begin{equation}\qquad
\bigl\|\mathbf X^{(\varepsilon)}(z)\mathbf x\bigr\|_2\le\bigl\|\mathbf
X^{(\varepsilon)}(z)\mathbf y\bigr\|_2+\bigl\|\mathbf X^{(\varepsilon)}(z)\bigr\|
\|\mathbf x-\mathbf y\|_2\le\gamma_1+\rho=\gamma_n/2.
\end{equation}
Hence,
%
%
\begin{equation}
\Pr(V)\le\exp\biggl\{-\frac{c_0}8 nq_n\biggr\}.
\end{equation}
Thus the lemma is proved.
\end{pf}
\begin{lem}\label{spread}Let $\delta, \rho\in(0,1)$. Let $\mathbf
x\in\mathit{Incomp}(\delta,\rho)$.
Then there exists a set $\sigma(\mathbf x)\subset\{1,\ldots,n\}$ of
cardinality $|\sigma(\mathbf x)|\ge\frac12 n\delta$ such that
%
%
\begin{equation}
\sum_{k\in\sigma(\mathbf x)}|x_k|^2\ge\frac12\rho^2
\end{equation}
and
%
%
\begin{equation}
\frac{\rho}{\sqrt{2n}}\le|x_k|\le\frac1{\sqrt{n\delta/2}}\qquad
\mbox{for any } k\in\sigma(\mathbf x),
\end{equation}
which we shall call ``spread set of $x$'' henceforth.
\end{lem}
\begin{pf}
See proof of Lemma 3.4 \cite{RV}, page 16. For the reader's
convenience we repeat this proof here.
Consider the subsets of $\{1,\ldots,n\}$ defined by
%
%
\begin{equation}
\sigma_1(\mathbf x):=\biggl\{k\dvtx|x_k|\le\frac{1}{\sqrt{\delta n/2}}\biggr\}
,\qquad \sigma_2(\mathbf x)=\biggl\{k\dvtx|x_k|\ge\frac{\rho}{\sqrt{2n}}\biggr\}
\end{equation}
and put $\sigma(\mathbf x)=\sigma_1(\mathbf x)\cap\sigma_2(\mathbf
x)$. Denote by $P_{\sigma(\mathbf x)}$ the orthogonal projection onto
$\mathbb R^{\sigma(\mathbf x)}$
in $\mathbb R^n$. By Chebyshev's inequality $|\sigma_1(\mathbf
x)^c|\le\delta n/2$.
Then $\mathbf y:={P}_{\sigma_1(\mathbf x)^c}\mathbf x\in
\mathit{Sparse}(\delta)$, so the incompressibility of $\mathbf x$
implies that
$\|{P}_{\sigma_1(\mathbf x)}\mathbf x\|_2=\|\mathbf x-\mathbf y\|
_2>\rho$. By the definition of $\sigma_2(\mathbf x)$, we have $\|
{P}_{\sigma_2(\mathbf x)^c}\mathbf x\|^2\le n\frac{\rho^2}{2n}=
\rho^2/2$.
Hence
%
%
\begin{equation}\label{in1}
\bigl\|{P}_{\sigma(\mathbf x)}\mathbf x\bigr\|_2^2\ge\bigl\|{P}_{\sigma_1(\mathbf
x)}\mathbf x\bigr\|_2^2-\bigl\|{P}_{\sigma_2(\mathbf x)}\mathbf x\bigr\|_2^2\ge\rho^2/2.
\end{equation}
Thus the lemma is proved.
\end{pf}
\begin{rem}\label{spre}If $\mathbf x\in\mathit{Incomp}(\delta\widehat
p_n,\rho)$ then there exists a set $\sigma(\mathbf x)$ with
cardinality $|\sigma(\mathbf x)|\ge\frac12 n\delta\widehat p_n$
such that
%
%
\begin{equation}
\frac{\rho}{\sqrt{2n}}\le|x_k|\le\frac1{\sqrt{n\delta\widehat p_n/2}}
\end{equation}
and
%
%
\begin{equation}
\bigl\|{P}_{\sigma(\mathbf x)}\mathbf x\bigr\|_2^2\ge\tfrac12\rho^2.
\end{equation}
\end{rem}

Let $Q(\eta)=\sup_{jk}\sup_{u\in\mathbb C}\Pr\{|X_{jk}-u|\le\eta
\}$.
Introduce the maximal concentration function of the weighed sums of the
rows of the matrix $(X_{jk})_{j,k=1}^n$,
%
%
\begin{equation}
p_\mathbf x(\eta)=\max_{j\in\{1,\ldots,n\}}\sup_{u\in\mathbb
C}\Pr\Biggl\{\Biggl|\sum_{k=1}^nX_{jk}\varepsilon_{jk}x_k-u\Biggr|\le\eta\Biggr\}.
\end{equation}
We shall now bound this concentration function and prove a
tensorization lemma for incompressible vectors.
\begin{lem}\label{new1}Let $\delta_n$ and $\rho_n$ be some functions
of $n$ such that $\rho_n,\delta_n\in(0,1)$.
Let $\eta_0$ and $r_0$ as in Lemma \ref{c01}.
Let $\mathbf x\in\mathit{Incomp}(\delta_n, \rho_n)$. Then there exists
positive constants $r_1$ and $r_2$ depending on $r_0$ such that for any
$0<\eta\le\eta_0$ we have
%
%
\begin{equation}
p_{\mathbf x}\bigl(\eta\rho_n/\sqrt{2n}\bigr)\le1- r_2\delta_n np_n
\end{equation}
for $n\delta_np_n\le1/3$ and
%
%
\begin{equation}
p_{\mathbf x}\bigl(\eta\rho_n/\sqrt{2n}\bigr)\le1-r_1<1
\end{equation}
for $n\delta_np_n> 1/3$.
\end{lem}
\begin{pf}Put $m=n\delta_n$. We have
%
%
\begin{eqnarray}
&&\sup_u\Pr\Biggl\{\Biggl|\sum_{k=1}^mX_{jk}\varepsilon_{jk}x_k-u\Biggr|\le\eta\rho
_n/\sqrt{2n}\Biggr\}\nonumber\\
&&\qquad\le\Pr\Biggl\{\sum_{k=1}^m\varepsilon_{jk}=0\Biggr\}\\
&&\qquad\quad{} +\Pr\Biggl\{\Biggl|\sum_{k=1}^mX_{jk}\varepsilon
_{jk}x_k-u\Biggr|\le\eta\rho_n/\sqrt{2n}; \sum_{k=1}^m\varepsilon
_{jk}\ge1\Biggr\}.\nonumber
\end{eqnarray}
Introduce $\sigma(\mathbf x):=\{k\in\{1,\ldots,n\}\dvtx\rho_n/{\sqrt
{2n}}\le|x_k|\le1/\sqrt{m/2}\}$.
Since $\mathbf x\in\break\mathit{Incomp}(\delta_n,\rho_n)$ the cardinality of
$\sigma(\mathbf x)$ is at least $m/2$.
Using that the concentration function of sum of independent random
variables is less then concentration function of its summands, we obtain
%
%
\begin{eqnarray}
&&\sup_u\Pr\Biggl\{\Biggl|\sum_{k=1}^mX_{jk}\varepsilon_{jk}x_k-u\Biggr|\le\eta\rho
_n/\sqrt{2n}\Biggr\}\nonumber\\[-8pt]\\[-8pt]
&&\qquad\le(1-p_n)^m+Q(\eta)\bigl(1-(1-p_n)^m\bigr).\nonumber
\end{eqnarray}
According to Lemma \ref{c01} in the \hyperref[app]{Appendix} for any
$\eta\le\eta
_0$, we have $Q(\eta)\le r_0<1$.
Assume that $mp_n\ge1/3$. Then we have
%
%
\begin{eqnarray}\qquad
\sup_u\Pr\Biggl\{\Biggl|\sum_{k=1}^mX_{jk}\varepsilon_{jk}x_k-u\Biggr|\le\eta\rho
_n/\sqrt{2n}\Biggr\}&\le& r_0+(1-r_0)e^{-mp_n}\nonumber\\
&\le&1-(1-e^{-1/3})(1-r_0)\\
&=&\!:1-r_1<1.\nonumber
\end{eqnarray}
If $mp_n\le1/3$ then $(1-p_n)^m\le1-mp_n/3$ and
%
%
\begin{eqnarray}
\sup_u\Pr\Biggl\{\Biggl|\sum_{k=1}^mX_{jk}\varepsilon_{jk}x_k-u\Biggr|\le\eta\rho
_n/\sqrt{2n}\Biggr\}&\le& 1-(1-r_0)mp_n/3\nonumber\\[-8pt]\\[-8pt]
&=&\!:1-r_2 mp_n.\nonumber
\end{eqnarray}
The lemma is proved.
\end{pf}

Now we state a tensorization lemma.
\begin{lem}\label{new2}
Let $\zeta_1,\ldots,\zeta_n$ be independent nonnegative random variables.
Assume that
%
%
\begin{equation}
\Pr\{\zeta_j\le\lambda_n\}\le1-q_n
\end{equation}
for some positive $q_n\in(0,1)$ and $\lambda_n>0$.
Then there exists positive absolute constants $K_1$ and $K_2$ such that
%
%
\begin{equation}
\Pr\Biggl\{\sum_{j=1}^n\zeta_j^2\le K_1^2nq_n\lambda_n^2\Biggr\}\le\exp\{
-K_2nq_n\}.
\end{equation}
\end{lem}
\begin{pf}We repeat the proof of Lemma 4.4 in \cite{LPRT}. Let
$t=K_1\sqrt{q_n}\lambda_n$. For any $\tau>0$ we have
%
%
\begin{equation}
\Pr\Biggl\{\sum_{j=1}^n\zeta_j^2\le nt^2\Biggr\}\le e^{n\tau}\prod
_{j=1}^n\mathbf{E}\exp\{-\tau\zeta_j^2/t^2\}.
\end{equation}
Furthermore,
%
%
\begin{eqnarray}
\mathbf{E}\exp\{-\tau\zeta_j^2/t^2\}&=&\int_0^{\infty}\Pr\bigl\{\exp\{
-\tau\zeta_j^2/t^2\}>s\bigr\}\,ds\nonumber\\
&=&
\int_0^{1}\Pr\bigl\{1/s>\exp\{\tau\zeta_j^2/t^2\} \bigr\}\,
ds\nonumber\\
&\le&\int_0^{\exp\{-\tau\lambda_n^2/t^2\}}ds+
\int_{\exp\{-\tau\lambda_n^2/t^2\}}^1(1-q_n)\,ds\\
&\le&1-q_n(1-\exp\{-\tau\lambda_n^2/t^2\})\nonumber\\
&=&1-q_n\bigl(1-\exp\{-\tau
/(K_1^2q_n)\}\bigr).\nonumber
\end{eqnarray}
Choosing $\tau:=q_n/4$ and $K_1^2:=\frac1{4\ln2}$, we get
%
%
\begin{equation}
\Pr\Biggl\{\sum_{j=1}^n\zeta_j^2\le nt^2\Biggr\}\le\exp\{-nq_n/2\}.
\end{equation}
Thus the lemma is proved.
\end{pf}

Recall that we assume $p_n^{-1} =O(n^{1-\theta}), 1\ge\theta>0$. For
this fixed
$\theta$ consider
$L:=[\frac1{\theta}]$. Hence by definition
$p_{n,l}:=(n\widehat p_n)^l p_n\to0, n \to\infty$ for
$l=1,\ldots,L-1$
and $\limsup_{n\to\infty}(np_n)^L p_n>0$. We put $p_{n,L}:=1$.

We shall assume that $n$ is large enough such that $(np_n)^Lp_n\ge
q_1>0$ for some constant $q_1>0$.
Starting with a decomposition of $\mathcal C_0:=\mathcal S^{(n-1)}$
into compressible vectors $\mathbf x$ in $\widehat{\mathcal C}_1:=
\mathcal C_0 \cap\mathit{Comp}(\delta_1
p_{n,1}, \rho_{n,1})$, where $p_{n,1}=\widehat p_n$, $\rho
_{n,1}=\gamma_0/(4K_n)$, and the constants $\gamma_0$ and $\delta_1$
are chosen as in Lemmas \ref{lem02} and \ref{lem031},
respectively. Then Lemma \ref{lem02} implies inequality (\ref{g01})
with $q_n$ replaced by $p_{n}$ and $\gamma_n$ replaced by $\gamma_0$.
Hence,\vspace*{1pt} using Lemma \ref{lem031}, one obtains the claim for the subset
of vectors $\widehat{\mathcal C}_1$.
The remaining vectors $\mathbf x$ in $\mathcal C_0$ lie in $\mathcal
C_1:= \mathit{Incomp}(\delta_1p_{n,1},\rho_{n,1})$. According to Lemmas
\ref{new1}, \ref{new2} inequality (\ref{g01}) holds again for these
vectors but with new parameters $q_n=np_n \delta_1 p_{n,1}$ and
$\gamma_n=c\rho_{n,1}\sqrt{\delta_1 p_{n,1}}$. Thus we may again subdivide
the vectors in $\mathcal C_1$ into the vectors within distance $\rho
_{n,2}$ from these
sparse ones, that is, $\widehat{\mathcal C}_2:= {\mathcal C}_1 \cap
\mathit{Comp}(\delta_2 p_{n,2}, \rho_{n,2})$ and the remaining ones, that
is, ${\mathcal C}_2:
= {\mathcal C}_1 \cap\mathit{Incomp}(\delta_2 p_{n,2},\rho_{n,2})$.
Iterating this procedure $L$ times we arrive at the incompressible set
${\mathcal C}_L$ of vectors $\mathbf x$ where Lemmas \ref{new1}, \ref
{new2} and Proposition~\ref{prop} yield the required bound of order
$\exp\{-\delta n\}$, for a sufficiently small absolute constant
$\delta>0$.

Summarizing, we will determine iteratively constants
$\delta_l, \rho_{n,l}$, for $l=1,\ldots,L$ and the following sets
of vectors:
%
%
\begin{equation}
\mathcal C_l:=\bigcap_{i=1}^{l}\mathit{Incomp}(\delta_i p_{n,i},\rho_{n,i})
\end{equation}
and
%
%
\begin{equation}
\widehat{\mathcal C}_l:=\mathcal C_{l-1}\cap\mathit{Comp}(\delta_l
p_{n,l},\rho_{n,l}) \qquad\mbox{with } \mathcal C_0=\mathcal S^{(n-1)}.
\end{equation}
Note that
%
%
\begin{equation}
\mathcal S^{(n-1)}=\bigcup_{l=1}^{L-1}\widehat{\mathcal C}_l\cup
\mathcal C_L.
\end{equation}
The main bounds to carry out this procedure are given in the following
Lemmas~\ref{spread1} and \ref{ziegel}.
\begin{lem}\label{spread1}
Let $\delta_{n},\rho_{n}\in(0,1)$ and let $\mathbf x\in\mathit
{Incomp}(\delta_{n},\rho_{n})$ and $\mathbf X^{(\varepsilon)}(z)$ be
a matrix as in Theorem \ref{thm1}.
Then there exist some positive constants $c_1$ and $c_2$ depending on
$K$, $r_0$, $\eta_0$ such
that for any $0<\tau\le\gamma_{n}$
%
%
\begin{equation}
\Pr\bigl\{\bigl\|\mathbf X^{(\varepsilon)}(z) \mathbf x\bigr\|_2\le\tau\bigr\}\le\exp
\bigl\{-c_1 n\bigl((p_n n\delta_{n})\wedge1\bigr)\bigr\}
\end{equation}
with
%
%
\begin{equation}\label{ga}
\gamma_{n}:=c_2{\rho_{n}\sqrt{\delta_{n}}},
\end{equation}
where $a\wedge b$ denotes the minimum of $a$ and $b$.
\end{lem}
\begin{pf} Assume at first that $n\delta_np_n\le1/3$.
According to Lemma \ref{new1}, we have, for any $j=1,\ldots,n$,
%
%
\begin{equation}
\sup_{u\in\mathbb C}\Pr\Biggl\{\Biggl|\sum_{k=1}^nX_{jk}\varepsilon
_{jk}x_k-u\Biggr|\le\eta_0\rho_n/\sqrt{2n}\Biggr\}\le1-r_1\delta_n np_n.
\end{equation}
Applying Lemma \ref{new2} with $q_n=r_1\delta_n np_n$, we get
%
%
\begin{equation}\qquad
\Pr\bigl\{\bigl\|\mathbf X^{(\varepsilon)}(z)\mathbf x\bigr\|_2\le\gamma_n/2
\mbox{ and } \bigl\|\mathbf X^{(\varepsilon)}(z)\bigr\|\le K_n\bigr\}\le\exp\{
-cn\delta_nnp_n\}.
\end{equation}

Consider now the case $n\delta_n p_n\ge1/3$. According to Lemma
\ref{new1}, we have
%
%
\begin{equation}
\sup_{u\in\mathbb C}\Pr\Biggl\{\Biggl|\sum_{k=1}^nX_{jk}\varepsilon
_{jk}x_k-u\Biggr|\le\eta_0\rho_n/\sqrt{2n}\Biggr\}\le1-r_1.
\end{equation}
Applying Lemma \ref{new2} with $q_n=r_1\delta_n np_n$, we get
%
%
\begin{equation}
\Pr\bigl\{\bigl\|\mathbf X^{(\varepsilon)}(z)\mathbf x\bigr\|_2\le\gamma_n/2
\mbox{ and } \bigl\|\mathbf X^{(\varepsilon)}(z)\bigr\|\le K_n\bigr\}\le\exp\{
-cn\}.
\end{equation}

This completes the proof of the lemma.
\end{pf}
\begin{lem}\label{ziegel} For $l=2,\ldots, L$ assume that $\delta_i,
\rho_{n,i}$ have been already determined for $i=1,\ldots,l-1$.
Then there exist absolute constants $\widehat c_l>0$ and $\overline
c_l>0$ and $\delta_l>0$ such that
%
%
\begin{eqnarray}
&&\Pr\Bigl\{\inf_{\mathbf x\in\widehat{\mathcal C_l}}\bigl\|\mathbf
X^{(\varepsilon)}(z) \mathbf x\bigr\|_2\le\gamma_{n,l} \mbox
{ and } \bigl\|\mathbf X^{(\varepsilon)}(z)\bigr\|\le
K_n\Bigr\}\nonumber\\[-8pt]\\[-8pt]
&&\qquad\le\exp\bigl\{-\overline c_l n\bigl(((n\widehat
p_n)^{l-1}p_n)\wedge1\bigr)\bigr\}\nonumber
\end{eqnarray}
with $\gamma_{n,l}$ defined by
%
%
\begin{equation}
\gamma_{n,l}=\widehat c_l \rho_{n,l-1}\sqrt{\delta_{l-1} p_{n,l-1}}
\end{equation}
and $\rho_{n,l}$ defined by
%
%
\begin{equation}
\rho_{n,l}:=\gamma_{n,l}/(4K_n),
\end{equation}
where $\widehat{\mathcal C}_l:=\mathcal C_{l-1}\cap\mathit{Comp}(\delta
_l p_{n,l},\rho_{n,l})$.
\end{lem}
\renewcommand{\theremark}{4.9}
\begin{remark}\label{rho} There exists some absolute constant $c>0$ that
%
%
\begin{equation}
\gamma_{n,L}\ge cn^{-L/2} \quad\mbox{and}\quad \rho_{n,L}\ge cn^{-(L+3)/2}.
\end{equation}
\end{remark}
\begin{pf}
Note that $p_{n,l}^{-1}=\mathcal O(n^{1-l\theta})$. This implies that
%
%
\begin{equation}\label{h01}
\gamma_{n,L}^{-1}=\rho_{n,1}^{-1}\mathcal O(n^{L-{L^2\theta}/2}).
\end{equation}
According to Lemmas \ref{lem02} and \ref{lem031}, we have $\rho
_{n1}^{-1}=\mathcal O(n^{({3-\theta})/2})$.
After simple calculations we get
%
%
\begin{equation}
\gamma_{n,L}^{-1}=\mathcal O(n^{L/2}).
\end{equation}
\upqed\end{pf}
\begin{pf*}{Proof of Lemma \ref{ziegel}}
To prove of this lemma we may use arguments similar to those in the
proofs of Lemmas 2.6 and 3.3 in \cite{RV}.
From $\mathbf x\in\mathcal C_l$ it follows that $\mathbf x\in\mathit
{Incomp}(\delta_{l-1}p_{n,l-1},\rho_{n,l-1})$.
Applying Lemma \ref{spread1} with $\delta_n=p_{n,l-1}$ and $\rho
_n=\rho_{n,l-1}$, we get
%
%
\begin{eqnarray}\label{k1}
&&\Pr\bigl\{\bigl\|\mathbf X^{(\varepsilon)}(z)\mathbf x\bigr\|_2\le\gamma
_{n,l} \mbox{ and } \bigl\|\mathbf X^{(\varepsilon)}(z)\bigr\|\le K_n\bigr\}
\nonumber\\[-8pt]\\[-8pt]
&&\qquad \le\exp\bigl\{-c_1n\bigl((np_n \widehat p_{n,l-1})\wedge
1\bigr)\bigr\}\nonumber
\end{eqnarray}
with
%
%
\begin{equation}
\gamma_{n,l}=c_2 \rho_{n,l-1}\sqrt{\delta_{l-1} p_{n,l-1}}.
\end{equation}
Inequality (\ref{k1}) and Lemma \ref{lem031} together imply
%
%
\begin{equation}\label{spread12}\qquad
\Pr\Bigl\{\inf_{\mathbf x\in{\mathcal C}_l}\bigl\|\mathbf X^{(\varepsilon
)}(z) \mathbf x\bigr\|_2\le\gamma_{n,l} \mbox{ and } \bigl\|\mathbf
X^{(\varepsilon)}(z)\bigr\|\le K_n\Bigr\}
\le\exp\{-c_1n \widehat p_{n,l}\}
\end{equation}
with $\delta_l$ defined in Lemma \ref{lem031} and
%
%
\begin{equation}
\rho_{n,l}:=\gamma_{n,l}/(4K_n).
\end{equation}

Thus the lemma is proved.
\end{pf*}

%
The next lemma gives an estimate of small ball probabilities adapted to
our case.
\begin{lem} \label{smalball}Let $\mathbf x\in\mathit{Incomp}(\delta
,\rho_{n,L})$.
Let $X_1,\ldots, X_n$ be random variables with zero mean and variance
at least 1.
Assume that the following condition holds:
%
%
\begin{equation}\label{unif}
L(M):=\max_{n\ge1}\max_{1\le k\le n}\mathbf{E}|X_k|^2I_{\{|X_k|>M\}
}\to0 \qquad\mbox{as }M\to\infty.
\end{equation}

Then there exist some constants $C>0$ depending on $\delta$ such that
for every
$\varepsilon>0$
%
%
\begin{equation}\label{eq:487}\qquad\quad
p_{\mathbf x}\bigl(\varepsilon\rho_{n,L}/\sqrt{2n}\bigr):=\sup_{v}\Pr\Biggl\{\Biggl|\sum
_{k=1}^nx_k\varepsilon_kX_k-v\Biggr|\le\varepsilon\rho_{n,L}/\sqrt{2n}\Biggr\}
\le\frac{C\sqrt{\ln n}}{\sqrt{np_n}}.
\end{equation}
\end{lem}
\begin{pf}
Put $L_1:=[-\log_2(\rho_{n,L}\sqrt{2\delta})]$. Note that
%
%
\begin{equation}
\frac{\rho_{n,L}}{\sqrt{2n}}\le\frac1{2^{L_1+1/2}\sqrt{n\delta
}}\le\frac{2\rho_{n,L}}{\sqrt{2n}}.
\end{equation}
According to Remark \ref{rho}, we have $\rho_{n,L}\ge cn^{-L/2}$.
This implies $L_1\le C\ln n$. Let $\sigma(\mathbf x)$ denote the
spread set of the vector $\mathbf x$, that is,
%
%
\begin{equation}
\sigma(\mathbf x):= \Biggl\{ k\dvtx\rho_{n,L}/\sqrt{2n}\le
|x_k|\le\sqrt{\frac2{n\delta}} \Biggr\}.
\end{equation}
By Lemma \ref{spread}, we have
%
%
\begin{equation}
|\sigma(\mathbf x)|\ge n\delta/2.
\end{equation}
We divide the spread interval of the vector $\mathbf x$ into $L_1+2$
intervals $\Delta_l$,
$l=0,\ldots,L_1+1$ by
%
%
\begin{eqnarray}
\Delta_0&:=&\biggl\{ k\dvtx\frac{\rho_{n,L}}{\sqrt{2n}}\le
|x_k|\le\frac1{2^{L_1+1/2}\sqrt{n\delta}} \biggr\},\\
\Delta_l&:=&\biggl\{ k\dvtx\frac{\sqrt{2}}{2^l\sqrt{n\delta}}\le
|x_k|\le\frac{\sqrt{2}}{2^{l-1}\sqrt{n\delta}} \biggr\}
,\qquad l=1,\ldots,L_1+1.
\end{eqnarray}
Note that there exists an $l_0\in\{0,\ldots,L_1+1\}$ such that
%
%
\begin{equation}
|\Delta_{l_0}|\ge n\delta/\bigl(2(L_1+2)\bigr)\ge Cn/\ln n.
\end{equation}
Let $\mathbf y=P_{\Delta_{l_0}}\mathbf x$. Put $a_{l}:={\min_{k\in
\Delta_{l}}}|x_k|$ and $b_{l}:=\max_{k\in\Delta_{l}}|x_k|$. Choose a
constant $M$ such that $L(M)\le1/2$. By the properties of
concentration functions, we have
%
%
\begin{equation}\label{eq:493}
p_{\mathbf x}\bigl(\varepsilon\rho_{n,L}/\sqrt{2n}\bigr)\le p_{\mathbf
y}\bigl(\varepsilon\rho_{n,L}/\sqrt{2n}\bigr)\le p_{\mathbf y}(Mb_{l_0}).
\end{equation}
By definition of $\Delta_{l_0}$, we have
%
%
\begin{equation}
\sum_{k\in\Delta_{l_0}}|x_k|^2\ge a_{l_0}^2|\Delta_{l_0}|\ge\rho
_{n,L}^2/(2n)|\Delta_{l_0}|
\end{equation}
and
%
%
\begin{equation}
\frac{a_{l_0}}{b_{l_0}}\ge\frac12.
\end{equation}
Define
%
%
\begin{equation}
D(\xi,\lambda)=\lambda^{-2}\mathbf{E}|\xi|^2I_{\{|\xi|<\lambda\}}
\end{equation}
and introduce for a random variable $\xi$, $\widetilde\xi:=\xi
-\widehat\xi$ where $\widehat\xi$ denotes an independent copy of
$\xi$. Put $\xi_k:=x_k\varepsilon_kX_k$.
We use the following inequality for a concentration function of a sum
of independent random variables:
%
%
\begin{equation}\label{petrov}
p_{\mathbf y}(Mb_{l_0})\le CMb_{l_0} \biggl(\sum_{k\in\Delta
_{l_0}}\lambda_k^2 D(\widetilde{\xi_k\varepsilon_k};\lambda
_k) \biggr)^{-1/2}
\end{equation}
with $\lambda_k\le Mb_{l_0}$.
See Petrov \cite{Petrov75}, page 43, Theorem 3. Put $\lambda_k=M|x_k|$.
It is straightforward to check that
%
%
\begin{equation}
\sum_{k\in\Delta_{l_0}}\lambda_k^2 D(\widetilde{\xi_k\varepsilon
_k};\lambda_k)\ge p_n \biggl(\sum_{k\in\Delta_{l_0}}|x_k|^2\bigl(\mathbf
{E}|X_k|^2-
L(M)\bigr) \biggr).
\end{equation}
This implies
%
%
\begin{equation}
\sum_{k\in\Delta_{l_0}}\lambda_k^2 D(\widetilde{\xi_k\varepsilon
_k};\lambda_k)\ge\frac{p_n}{2}\sum_{k\in\Delta_{l_0}}|x_k|^2
\ge\frac{p_n}{2}|\Delta_{l_0}|a_{l_0}^2.
\end{equation}
Combining this inequality with (\ref{petrov}) and (\ref{eq:493}) we obtain
%
%
\begin{equation}
p_{\mathbf x}\bigl(\varepsilon\rho_{n,L}/\sqrt{2n}\bigr)\le\frac
{CMb_{l_0}}{\sqrt{|\Delta_{l_0}|p_n}a_{l_0}}\le\frac{CM}{\sqrt
{|\Delta_{l_0}|p_n}}\le\frac{C\sqrt{\ln n}}{\sqrt{np_n}}.
\end{equation}
The last relation concludes the proof.
\end{pf}

\subsection*{Invertibility for the incompressible vectors via distance}
\begin{lem}\label{invert}
Let $\mathbf X_1,\mathbf X_2,\ldots,\mathbf X_n$ denote the columns of
$\sqrt{np_n}\mathbf X^{(\varepsilon)}(z)$, and let $\mathcal H_k$
denotes the span of
all column vectors except the $k$th.
Then for every $\delta, \rho\in(0,1)$ and every $\eta>0$ one has
\begin{eqnarray*}
&&\Pr\Bigl\{\inf_{\mathbf x\in{\mathcal C}_L}\bigl\|\mathbf
X^{(\varepsilon)}(z)
\mathbf x\bigr\|_2<\eta\bigl(\rho
_{n,L}/\sqrt{n}\bigr)^2/\sqrt{np_n} \Bigr\}\\
&&\qquad
\le\frac{1}{n\delta_L }\sum_{k=1}^n\Pr\bigl\{\operatorname{dist}(\mathbf
X_k,\mathcal H_k)<\eta\rho_{n,L}/\sqrt{n}\bigr\}.
\end{eqnarray*}
\end{lem}
\begin{pf}Note that
%
%
\begin{eqnarray}\label{in00}\quad
&&\Pr\Bigl\{\inf_{\mathbf x\in\widehat{\mathcal C}_L}\bigl\|\mathbf
X^{(\varepsilon)}(z) \mathbf x\bigr\|_2<\eta\bigl(\rho_{n,L}/\sqrt
{n}\bigr)^2/\sqrt{np_n} \Bigr\}\nonumber\\[-8pt]\\[-8pt]
&&\qquad \le
\Pr\Bigl\{\inf_{\mathbf x\in\mathit{Incomp}(\delta_L,\rho_{n,L})}\bigl\|
\mathbf X^{(\varepsilon)}(z) \mathbf x\bigr\|_2<\eta\bigl(\rho_{n,L}/\sqrt
{n}\bigr)^2/\sqrt{np_n} \Bigr\}.\nonumber
\end{eqnarray}
For the upper bound of the r.h.s. of (\ref{in00}) (see \cite{RV},
proof of Lemma 3.5). For the reader's convenience we repeat this proof.
Introduce the matrix $\mathbf G:=\sqrt{np_n}\mathbf X^{(\varepsilon
)}(z)$. Recall that $\mathbf X_1,\ldots,\mathbf X_n$ denote the column
vector of the matrix $\mathbf G$ and $\mathcal H_k$ denotes the span of
all column vectors except the $k$th. Writing $\mathbf G\mathbf x=\sum
_{k=1}^nx_k\mathbf X_k$,
we have
%
%
\begin{equation}\label{dist}
\|\mathbf G\mathbf x\|\ge\max_{k=1,\ldots,n}\operatorname
{dist}(x_k\mathbf
X_k,\mathcal H_k)={\max_{k=1,\ldots,n}}|x_k|\operatorname{dist}(\mathbf
X_k,\mathcal H_k).
\end{equation}
Put
%
%
\begin{equation}
p_k:=\Pr\bigl\{\operatorname{dist}(\mathbf X_k,\mathcal H_k)<\eta\rho
_{n,L}/\sqrt{n} \bigr\}.
\end{equation}
Then
%
%
\begin{equation}
\mathbf{E} \bigl|\bigl\{k\dvtx\operatorname{dist}(\mathbf X_k,\mathcal H_k)<\eta
\rho
_{n,L}/\sqrt{n}\bigr\} \bigr|=\sum_{k=1}^np_k.
\end{equation}
Denote by $U$ the event that the set $\sigma_1:=\{k\dvtx\operatorname
{dist}(\mathbf
X_k,H_k)\ge\eta\rho_{n,L}/\sqrt{n}\}$ contains more than $(1-\delta
_L)n$ elements. Then by Chebyshev's inequality
%
%
\begin{equation}
\Pr\{U^c\}\le\frac1{n\delta_L}\sum_{k=1}^np_k.
\end{equation}
On the other hand, for every incompressible vector $\mathbf x$, the set
$\sigma_2(\mathbf x):=\{k\dvtx|x_k|\ge\rho_{n,L}/\sqrt n\}$ contains
at least $n\delta_L$ elements. (Otherwise, since\break $\|P_{\sigma
_2(\mathbf x)^c}\mathbf x\|_2\le\rho_{n,L}$, we have $\|\mathbf
x-\mathbf y\|_2\le\rho_{n,L}$ for the sparse vector $\mathbf
y:=P_{\sigma_2(\mathbf x)}\mathbf x$, which would contradict the
incompressibility of $\mathbf x$.)

Assume that the event $U$ occurs. Fix any incompressible vector
$\mathbf x$. Then $|\sigma_1|+|\sigma_2(\mathbf x)|>(1-\delta
_{L})n+n\delta_L>n$, so the sets $\sigma_1$ and $\sigma_2(\mathbf
x)$ have nonempty intersection. Let $k\in\sigma_1\cap\sigma
_2(\mathbf x)$. Then by (\ref{dist}) and by definitions of the sets
$\sigma_1$ and $\sigma_2(\mathbf x)$, we have
%
%
\begin{equation}
\|\mathbf G\mathbf x\|_2\ge|x_k|\operatorname{dist}(\mathbf X_k,\mathcal
H_k)\ge\eta(\rho_{n,L}n^{-1/2})^2.
\end{equation}
Summarizing we have shown that
%
%
\begin{equation}\qquad\quad
\Pr\Bigl\{\inf_{\mathbf x\in\mathit{Incomp}(\delta_L,\rho_{n,L})}\|\mathbf
G\mathbf x\|_2\le\eta(\rho_{n,L}n^{-1/2})^2\Bigr\}\le\Pr\{U^c\}\le
\frac1{n\delta_L}\sum_{k=1}^np_k.
\end{equation}

This completes the proof.
\end{pf}

We now reformulate Lemma 3.6 from \cite{RV}. Let $\mathbf X_n^*$ be
any unit vector orthogonal to $\mathbf X_1,\ldots, \mathbf X_{n-1}$.
Consider the subspace $\mathcal H_n=\operatorname{span}(\mathbf
X_1,\ldots,
\mathbf X_{n-1})$.
\begin{lem}\label{normal} Let $\delta_l, \rho_l, c_l$, $l=1,\ldots,
L-1$, be as in Lemma \ref{lem031} and $\delta_L$, $\rho_L,\overline c_L$ as
in Lemma \ref{ziegel}. Then there exists an absolute constant
$\widehat c_L>0$ such that
%
%
\begin{equation}
\Pr\bigl\{\mathbf X^*\notin{\mathcal C}_L\mbox{ and } \bigl\|\mathbf
X^{(\varepsilon)}(z)\bigr\|\le K_n \bigr\}\le\exp\{-\widehat c_Lnp_n\}.
\end{equation}
\end{lem}
\begin{pf}
Note that
%
%
\begin{equation}
\mathcal S^{(n-1)}=\bigcup_{l=1}^{L-1}\widehat{\mathcal C}_l\cup
{\mathcal C}_L.
\end{equation}
The event $\{\mathbf X^*\notin{\mathcal C}_L \mbox{ and } \|
\mathbf X^{(\varepsilon)}(z)\|\le K_n\}$ implies that the event
%
%
\begin{equation}\quad
\mathcal E:=\Bigl\{\inf_{\mathbf x \in\bigcup_{l=1}^{L-1}\widehat{\mathcal
C}_l\dvtx\|\mathbf x\|_2=1}\bigl\|\mathbf X^{(\varepsilon)}(z)\mathbf x\bigr\|
_2\le c \mbox{ and }
\bigl\|\mathbf X^{(\varepsilon)}(z)\bigr\|\le K_n\Bigr\}
\end{equation}
occurs for any positive $c$. This implies, for $c>0$,
%
%
\begin{eqnarray}\qquad
&&\Pr\bigl\{\mathbf X^*\notin{\mathcal C}_L \mbox{ and } \bigl\|\mathbf
X^{(\varepsilon)}(z)\bigr\|\le K_n\bigr\}\\
&&\qquad\le\sum_{l=1}^{L-1}\Pr\Bigl\{\inf_{\mathbf x \in\widehat{\mathcal
C}_l\dvtx\|\mathbf x\|_2=1}\bigl\|\mathbf X^{(\varepsilon)}(z)\mathbf x\bigr\|\le
c \mbox{ and } \bigl\|\mathbf X^{(\varepsilon)}(z)\bigr\|\le K_n\Bigr\}.
\end{eqnarray}
Now choose $c:=\min\{\gamma_{n,l}, l=1,\ldots, L-1\}$. Applying
Lemma \ref{ziegel}
proves the claim.
\end{pf}
\begin{lem}\label{dist1}Let $\mathbf X^{(\varepsilon)}(z)$ be a
random matrix as in Theorem \ref{sparse}.
Let $\mathbf X_1,\ldots,\mathbf X_n$ denote column vectors of the
matrix $\sqrt{np_n}\mathbf X^{(\varepsilon)}(z)$, and consider the
subspace $\mathcal H_n=\operatorname{span}(\mathbf X_1,\ldots,\mathbf X_{n-1})$.
Let $K_n=Kn\sqrt{p_n}$. Then we have
%
%
\begin{equation}\quad 
\Pr\bigl\{\operatorname{dist}(\mathbf X_n,\mathcal H_n)<\rho_{n,L}/\sqrt n
\mbox{ and } \bigl\|\mathbf X^{(\varepsilon)}(z)\bigr\|\le K_n\bigr\}\le\frac
{C\sqrt{\ln n}}{\sqrt{np_n}}.
\end{equation}
\end{lem}
\begin{pf}We repeat Rudelson and Vershynin's proof of Lemma 3.8 in
\cite{RV}.
Let $\mathbf X^*$ be any unit vector orthogonal to $\mathbf X_1,\mathbf
X_2,\ldots,\mathbf X_{n-1}$. We can choose $\mathbf X^*$ so that it is
a random vector that depends on $\mathbf X_1,\mathbf X_2,\ldots
,\mathbf X_{n-1}$ only and is independent of~$\mathbf X_n$. We have
\[
\operatorname{dist}(\mathbf X_n,\mathcal H_n)\ge|\langle\mathbf X_n,\mathbf X^*\rangle|.
\]
We denote the probability with respect to $\mathbf X_n$ by ${\Pr}_n$
and the expectation with respect to $\mathbf X_1,\ldots,\mathbf
X_{n-1}$ by
$\mathbf{E}_{1,\ldots,n-1}$. Then
%
%
\begin{eqnarray}
&&\Pr\bigl\{\operatorname{dist}(\mathbf X_n,\mathcal H_n)<\rho_{n,L}/\sqrt n
\mbox{ and } \bigl\|\mathbf X^{(\varepsilon)}(z)\bigr\|
\le K_n
\bigr\}\nonumber\\
&&\qquad\le
\mathbf{E}_{1,\ldots,n-1}{\Pr}_n\bigl\{|\langle\mathbf X^*,\mathbf X_n\rangle|\le
\rho_{n,L}/\sqrt n \mbox{ and } \mathbf X^*\in{\mathcal
C}_L\bigr\}
\\
&&\qquad\quad{} +\Pr\bigl\{\mathbf X^*\notin{\mathcal
C}_L \mbox{ and } \bigl\|\mathbf X^{(\varepsilon)}(z)\bigr\|\le K_n\bigr\}.\nonumber
\end{eqnarray}
According to Lemma \ref{normal}, the second term in the right-hand
side of the last inequality is less then $\exp\{-\widehat c_Ln\}$.
Since the vectors $\mathbf X^*=(a_1,\ldots,a_n)\in\mathcal S^{(n-1)}$
and $\mathbf X_n=(\varepsilon_1\xi_1,\ldots,\varepsilon_n\xi_n)$
are independent,
we may use small ball probability estimates. We have
\[
S=\langle\mathbf X_n,\mathbf X^*\rangle=\sum_{k=1}^na_k\varepsilon_k\xi_k.
\]
Let $\sigma$ denote the spread set of $\mathbf X^*$ as in Lemma \ref
{spread}. Let $P_{\sigma}$ denote the orthogonal projection onto
$\mathbb R^{\sigma}$ in $\mathbb R^n$. Denote by $S_{\sigma}=\sum
_{k\in\sigma}\varepsilon_k a_k\xi_k$. Using the properties of
concentration functions, we get
\begin{eqnarray*}
{\Pr}_n\bigl\{|\langle\mathbf X_n,\mathbf X^*\rangle|\le\rho_{n,L}/\sqrt
n\bigr\}
&\le&\sup_v{\Pr}_n\bigl\{|S-v|\le\rho_{n,L}/\sqrt n\bigr\}\\
&\le&\sup_v{\Pr}_n\bigl\{|S_{\sigma}-v|\le\rho_{n,L}/\sqrt n\bigr\}.
\end{eqnarray*}
By Lemma \ref{smalball}, we have for some absolute constant $C>0$
%
%
\begin{equation}\label{101}
{\Pr}_n\bigl\{|\langle\mathbf X_n,\mathbf X^*\rangle |\le\rho_{n,L}/\sqrt n\bigr\}\le\frac
{C\sqrt{\ln n}}{\sqrt{np_n}}.
\end{equation}

Thus the lemma is proved.
\end{pf}
\begin{lem}\label{comp}Let $\mathbf X^{(\varepsilon)}(z)$ be a random
matrix as in Theorem \ref{thm1}.
Let $\delta_L, \rho_{n,L}\in(0,1)$.
Let $\mathbf X_1,\ldots,\mathbf X_n$ denote column vectors of matrix
$\sqrt{np_n}\mathbf X^{(\varepsilon)}(z)$.
Let $K_n=Kn\sqrt{p_n}$ with $K\ge1$. Then we have
\[
\Pr\Bigl\{\inf_{\mathbf x\in{\mathcal C}_L}\bigl\|\mathbf X^{(\varepsilon
)}(z)\mathbf x\bigr\|_2<\rho_{n,L}^2 /n\Bigr\}
\le\Pr\bigl\{\bigl\|\mathbf X^{(\varepsilon)}(z)\bigr\|>K_n\bigr\}+\frac{C\sqrt{\ln
n}}{\sqrt{np_n}}.
\]
\end{lem}
\begin{pf}Note that
%
%
\begin{eqnarray}
&&\Pr\Bigl\{\inf_{\mathbf x\in{\mathcal C}_L}\bigl\|\mathbf X^{(\varepsilon
)}(z) \mathbf x\bigr\|_2<\rho_{n,L}^2 /n\Bigr\}
\nonumber\\
&&\qquad\le
\Pr\Bigl\{\inf_{\mathbf x\in{\mathcal C}_L}\bigl\|\mathbf X^{(\varepsilon
)}(z)\mathbf x\bigr\|_2<\rho_{n,L}^2/ n
\mbox{ and }
\bigl\|\mathbf X^{(\varepsilon)}(z)\bigr\|\le K_n\Bigr\}\\
&&\qquad\quad{} +
\Pr\bigl\{\bigl\|\mathbf X^{(\varepsilon)}(z)\bigr\|>K_n\bigr\}.\nonumber
\end{eqnarray}
Applying Lemma \ref{invert} with $\eta=\sqrt{p_n}$, we get
\[
\Pr\biggl\{\inf_{\mathbf x\in{\mathcal C}_L}\bigl\|\mathbf X^{(\varepsilon
)}(z) \mathbf x\bigr\|_2<\frac{\rho_{n,L}^2}{n}\biggr\}
\le
\frac{1}{n \delta_L }\sum_{k=1}^n\Pr\biggl\{\operatorname{dist}(\mathbf
X_k,\mathcal H_k)<\frac{\rho_{n,L}\sqrt{p_n}}{\sqrt{n}}\biggr\}.
\]
Applying Lemma \ref{dist1}, we obtain
%
%
\begin{equation}
\Pr\Bigl\{\inf_{\mathbf x\in{\mathcal C}_L}\bigl\|\mathbf X^{(\varepsilon
)}(z) \mathbf x\bigr\|_2<\rho_{n,L}^2 / n\Bigr\}\le\frac{C\sqrt{\ln
n}}{\sqrt{n p_n}}.
\end{equation}
Thus the lemma is proved.
\end{pf}
%
\begin{pf*}{Proof of Theorem \ref{thm1}}
By definition of the minimal singular value, we have
\begin{eqnarray*}
&&\Pr\bigl\{s_n^{(\varepsilon)}(z)\le\rho_{n,L}^2 /{n} \mbox
{ and } s_1^{(\varepsilon)}(z)\le
K_n\bigr\}\\
&&\qquad \le
\Pr\bigl\{\mbox{there exists } \mathbf x\in\mathcal S^{(n-1)}\dvtx\bigl\|\mathbf
X^{(\varepsilon)}(z)\mathbf x\bigr\|_2 \le\rho_{n,L}^2 /{n}
\mbox{ and } s_1^{(\varepsilon)}(z)\le K_n\bigr\}.
\end{eqnarray*}
Furthermore, using the decomposition of the sphere $\mathcal
S^{(n-1)}=\bigcup_{l=1}^{L-1}\widehat{\mathcal C}_l\cup\mathcal C_L$
into compressible and incompressible vectors, we get
%
%
\begin{eqnarray}
&&\Pr\bigl\{s_n^{(\varepsilon)}(z)\le\rho_{n,L}^2 /{n}\mbox{ and }
s_1^{(\varepsilon)}(z)\le K_n\bigr\}
\nonumber\\
&&\qquad\le\sum_{l=1}^{L-1}
\Pr\Bigl\{\inf_{\mathbf x\in\widehat{\mathcal C}_l}\bigl\|\mathbf
X^{(\varepsilon)}(z)\mathbf x\bigr\|_2 \le\rho_{n,L}^2 /{n}
\mbox{ and } s_1^{(\varepsilon)}(z)\le K_n\Bigr\}\\
&&\qquad\quad{} +\Pr\Bigl\{\inf_{\mathbf x\in{\mathcal C}_L}\bigl\|\mathbf X^{(\varepsilon
)}(z)\mathbf x\bigr\|_2 \le\rho_{n,L}^2 /{n}\mbox{ and }
s_1^{(\varepsilon)}(z)\le K_n\Bigr\}.\nonumber
\end{eqnarray}
According to Lemma \ref{ziegel}, we have
\[
\Pr\Bigl\{\inf_{\mathbf x\in\widehat{\mathcal C}_l}\bigl\|\mathbf
X^{(\varepsilon)}(z)\mathbf x\bigr\|_2 \le\rho_{n,L}^2
/{n}\mbox{ and }
s_1^{(\varepsilon)}(z)\le K_n\Bigr\}\le\exp\{-c_l np_n(n\widehat
p_n)^{l-1}\}.
\]
Lemmas \ref{comp} and \ref{ziegel} together imply that
%
%
\begin{eqnarray}\label{eq:4120}\qquad
&&\Pr\Bigl\{\inf_{\mathbf x\in{\mathcal C}_L}\bigl\|\mathbf X^{(\varepsilon
)}(z)\mathbf x\bigr\|_2 \le\rho_{n,L}^2 /{n}\mbox{ and }
s_1^{(\varepsilon)}(z)\le K_n\Bigr\}
\nonumber\\
&&\qquad
\le
\Pr\Bigl\{\inf_{\mathbf x\in\mathit{Incomp}(\delta_L,\rho_{n,L})}\bigl\|\mathbf
X^{(\varepsilon)}(z)\mathbf x\bigr\|_2 \le\rho_{n,L}^2 /{n}
\mbox{ and } s_1^{(\varepsilon)}(z)\le K_n\Bigr\}
\\
&&\qquad
\le
\frac{C\sqrt{\ln n}}{\sqrt{np_n}}+\exp\{-\widehat c_Ln\}.\nonumber
\end{eqnarray}
The last two inequalities together imply the result.
\end{pf*}
\begin{rem}\label{rem:relax}
To relax the condition $p_n^{-1}=\mathcal O(n^{1-\theta})$ of Theorem
\ref{thm1} to $p_n^{-1}= o(n/\ln^2{n})$ we should put $L=\ln{n}$.
Then the value $L_1$ in Lemma \ref{smalball} is at most $C(\ln
{n})^2$, and hence we get the bound $C\ln{n}/\sqrt{np_n}$ in (\ref
{eq:487}). This yields the bound $C\ln{n}/\sqrt{np_n}+\exp\{
-\widehat c_Ln\}$ in (\ref{eq:4120}). Thus Theorem \ref{thm1} holds
with $B$ chosen to be of order $C\ln{n}$.
\end{rem}

\section{Proof of the main theorem}\label{proof}
In this section we give the proof of Theorem~\ref{sparse}. Theorem
\ref{thm0} follows from Theorem \ref{sparse} with $p_n=1$. Let
$\gamma:=\frac1{3}$ and let $R>0$ and $k_1$ be defined as in Lemma
\ref{compact} with $q=18$.
Using the notation of Theorem \ref{thm1} we introduce for any $z\in
\mathbb C$ and absolute constant $c>0$ the set
$\Omega_n(z)=\{\omega\in\Omega\dvtx c/n^B\le s_n^{(\varepsilon)}(z),
s_1{(\varepsilon)}\le n\sqrt{p_n},
|\lambda_{k_1}^{(\varepsilon)}|\le R\}$.
According to Lemma \ref{largeval}
\[
\Pr\bigl\{s_1^{(\varepsilon)}(\mathbf X)\ge n\sqrt{p_n}\bigr\}\le C(np_n)^{-1}.
\]
According to Theorem \ref{thm1} with $\varepsilon=c$, we have
\[
\Pr\bigl\{c/n^B\ge s_n^{(\varepsilon)}(z)\bigr\}\le\frac{C\sqrt{\ln
n}}{\sqrt{np_n}}+\Pr\bigl\{s_1^{(\varepsilon)}\ge n\sqrt{p_n}\bigr\}.
\]
According to Lemma \ref{compact} with $q=18$, we have
%
%
\begin{equation}
\Pr\bigl\{\bigl|\lambda_{k_1}^{(\varepsilon)}\bigr|\le R\bigr\}\le C\Delta_n^{\gamma
}\le C\bigl[\varphi\bigl(\sqrt{np_n}\bigr)\bigr]^{-{1}/{18}}.
\end{equation}

These inequalities imply
%
%
\begin{equation}\label{truncation}
\Pr\{\Omega_n(z)^{c}\}\le\bigl(\varphi\bigl(\sqrt{np_n}\bigr)\bigr)^{-1/{18}}.
\end{equation}
Let $r=r(n)$ be such that $r(n)\to0$ as $n\to\infty$. A more
specific choice will be made later.
Consider the potential $U_{\mu_n}^{(r)}$.
We have
\begin{eqnarray*}
U_{\mu_n}^{(r)}&=&-\frac1n\mathbf{E}\log\bigl|\det\bigl(\mathbf
X^{(\varepsilon)}-z\mathbf I-r\xi\mathbf I\bigr)\bigr|\\
&=&-\frac1n\sum_{j=1}^n\mathbf{E}\log\bigl|\lambda_j^{(\varepsilon
)}-r\xi-z\bigr|I_{\Omega_n(z)}\\
&&{}-\frac1n\sum_{j=1}^n\mathbf{E}\log\bigl|\lambda_j^{(\varepsilon)}-r\xi
-z\bigr|I_{\Omega_n^{(c)}(z)}\\
&=&
\overline U{}^{(r)}_{\mu_n}+\widehat U_{\mu_n}^{(r)},
\end{eqnarray*}
where
$I_A$ denotes an indicator function of an event $A$ and ${\Omega
_n(z)}^{c}$ denotes the
complement of $\Omega_n(z)$.
\begin{lem}\label{lem5.1}Assuming the conditions of
Theorem \ref{thm1}, for $r$ such that
\[
\ln(1/r) \bigl(\varphi\bigl(\sqrt{np_n}\bigr)\bigr)^{-1/19}\to\infty\qquad\mbox{as }
n\to\infty
\]
we have
%
%
\begin{equation}\label{0*}
\widehat U_{\mu_n}^{(r)}\to0\qquad\mbox{as }n\to\infty.
\end{equation}
\end{lem}
\begin{pf}
By definition, we have
%
%
\begin{equation}\label{1*}
\widehat U_{\mu_n}^{(r)}=-\frac1n\sum_{j=1}^n\mathbf{E}\log
\bigl|\lambda_j^{(\varepsilon)}-r\xi-z\bigr|I_{\Omega_n^{(c)}(z)}.
\end{equation}
Applying Cauchy's inequality, we get, for any $\tau>0$,
%
%
\begin{eqnarray}\label{2*}\quad
\bigl|\widehat U_{\mu_n}^{(r)}\bigr|&\le&\frac1n\sum_{j=1}^n\mathbf{E}^{
1/({1+\tau})}\bigl|\log\bigl|\lambda_j^{(\varepsilon)}-r\xi-z\bigr|\bigr|^{1+\tau}
\bigl(\Pr\bigl\{\Omega_n^{(c)}\bigr\}\bigr)^{{\tau}/({1+\tau
})}\nonumber\\[-8pt]\\[-8pt]
&\le&
\Biggl(\frac1n\sum_{j=1}^n\mathbf{E}\bigl|\log\bigl|\lambda_j^{(\varepsilon
)}-r\xi-z\bigr|\bigr|^{1+\tau} \Biggr)^{1/({1+\tau})} \bigl(\Pr\bigl\{\Omega
_n^{(c)}\bigr\} \bigr)^{{\tau}/({1+\tau})}.\nonumber
\end{eqnarray}
Furthermore, since $\xi$ is uniformly distributed in the unit disc and
independent of~$\lambda_j$, we may write
\begin{eqnarray*}
\mathbf{E} \bigl|{\log}|\lambda_j-r\xi-z| \bigr|^{1+\tau}&=&\frac
1{2\pi}\mathbf{E}\int_{|\zeta|\le1} \bigl|\log\bigl|\lambda
_j^{(\varepsilon)}
-r\zeta-z\bigr| \bigr|^{1+\tau}\,d\zeta
\\
&=&
\mathbf{E}J_1^{(j)}+\mathbf{E}J_2^{(j)}+\mathbf{E}J_3^{(j)},
\end{eqnarray*}
where
\begin{eqnarray*}
J_1^{(j)}&=&\frac1{2\pi}\int_{|\zeta|\le1, |\lambda
_j^{(\varepsilon)}-r\zeta-z|\le\varepsilon}\bigl|\log\bigl|\lambda
_j^{(\varepsilon)}-r\zeta-z\bigr|\bigr|^{1+\tau}\,d\zeta,\\
J_2^{(j)}&=&\frac1{2\pi}\int_{|\zeta|\le1, 1/{\varepsilon
}>|\lambda_j^{(\varepsilon)}-r\zeta-z|>\varepsilon}\bigl|\log\bigl|\lambda
_j^{(\varepsilon)}-r\zeta-z\bigr|\bigr|^{1+\tau}\,d\zeta,\\
J_3^{(j)}&=&\frac1{2\pi}\int_{|\zeta|\le1, |\lambda_j-r\zeta
-z|>1/{\varepsilon}}\bigl|\log\bigl|\lambda_j^{(\varepsilon)}-r\zeta
-z\bigr|\bigr|^{1+\tau}\,d\zeta.
\end{eqnarray*}
Note that
\[
\bigl|J_2^{(j)}\bigr|\le\log\biggl(\frac1{\varepsilon} \biggr).
\]
Since for any $b>0$, the function $-u^b\log u$ is not decreasing on the
interval $[0,\exp\{-\frac1{b}\}]$,
we have for $0<u\le\varepsilon<\exp\{-\frac1{b}\}$,
\[
-\log u\le\varepsilon^{b}u^{-b}\log\biggl(\frac1{\varepsilon
} \biggr).
\]
Using this inequality, we obtain, for $b(1+\tau)<2$,
%
%
\begin{eqnarray}\label{finish1}
\bigl|J_1^{(j)}\bigr|&\le&\frac1{2\pi}\varepsilon^{b(1+\tau)} \biggl(\log
\biggl(\frac1{\varepsilon} \biggr) \biggr)^{1+\tau}\nonumber\\[-8pt]\\[-8pt]
&&{}\times
\int_{|\zeta|\le1, |\lambda_j^{(\varepsilon)}-r\zeta-z|\le
\varepsilon}\bigl|\lambda_j^{(\varepsilon)}-r\zeta-z\bigr|^{-b(1+\tau
)}\,d\zeta\nonumber\\
&\le&\frac1{2\pi r^2}\varepsilon^{b}\log\biggl(\frac1{\varepsilon
}
\biggr)\int_{|\zeta|\le\varepsilon}|\zeta|^{-b(1+\tau)}\,d\zeta\nonumber\\[-8pt]\\[-8pt]
&\le& C(\tau, b)\varepsilon^{2}r^{-2} \biggl(\log\biggl(\frac
1{\varepsilon} \biggr) \biggr)^{1+\tau}.\nonumber
\end{eqnarray}
If we choose $\varepsilon=r$, then we get
%
%
\begin{equation}\label{finish2}
\bigl|J_1^{(j)}\bigr|\le C(\tau, b) \biggl(\log\biggl(\frac1{r} \biggr)
\biggr)^{1+\tau}.
\end{equation}
The following bound holds for $\frac1n\sum_{j=1}^n\mathbf
{E}J_3^{(j)}$. Note that
$|{\log x}|^{1+\tau}\le\varepsilon^2\times\break|{\log\varepsilon}
|^{1+\tau} x^2$ for $x\ge\frac1{\varepsilon}$ and
sufficiently small $\varepsilon$.
Using this inequality, we obtain
%
%
\begin{eqnarray}\label{finish3}\qquad
\frac1n\sum_{j=1}^n\mathbf{E}J_3^{(j)}&\le& C(\tau)\varepsilon
^2|{\log\varepsilon}|\frac1n\sum_{j=1}^n \mathbf{E}\bigl|\lambda
_j^{(\varepsilon)}-r\zeta-z\bigr|^2\nonumber\\
&\le& C(\tau)(1+|z|^2+r^2)\varepsilon^2|{\log\varepsilon}|
\\
&\le& C(\tau)(2+|z|^2)r^2|{\log r}|.\nonumber
\end{eqnarray}

Inequalities (\ref{finish1})--(\ref{finish3}) together imply that
%
%
\begin{equation}\label{3*}
\Biggl|\frac1n\sum_{j=1}^n\mathbf{E}\bigl|{\log}\bigl|\lambda_j^{(\varepsilon
)}-r\xi-z\bigr|\bigr|^{1+\tau}\Biggr|\le
C \biggl(\log\biggl(\frac1{r} \biggr) \biggr)^{1+\tau}.
\end{equation}
Furthermore, inequalities (\ref{truncation}), (\ref{1*}), (\ref{2*})
and (\ref{3*}) together imply
\[
\bigl|\widehat U_{\mu_n}^{(r)}\bigr|\le C \biggl(\log\biggl(\frac1{r}
\biggr) \biggr) \bigl(C\bigl(\varphi\bigl(\sqrt{np_n}\bigr)\bigr)^{-1/{18}}
\bigr)^{{\tau}/({1+\tau})}.
\]
We choose $\tau=18$ and rewrite the last inequality as follows:
\[
\bigl|\widehat U_{\mu_n}^{(r)}\bigr|\le C \biggl(\log\biggl(\frac1{r}
\biggr) \biggr)\bigl(\varphi\bigl(\sqrt{np_n}\bigr)\bigr)^{-{1}/{19}}
\le C \biggl(\log\biggl(\frac1{r} \biggr) \biggr)\bigl(\varphi\bigl(\sqrt
{np_n}\bigr)\bigr)^{-{1}/{19}}.
\]
If we choose $r=\frac1{\sqrt{np_n}}$ we obtain
$\log(1/r)((\varphi(\sqrt{np_n}))^{-{1}/{19}}\to0$, then (\ref
{0*}) holds and the
lemma is proved.
\end{pf}

We shall investigate $\overline U{}^{(r)}_{\mu_n}$ now.
We may write
%
%
\begin{eqnarray}
\overline U{}^{(r)}_{\mu_n}&=&-\frac1n\sum_{j=1}^n\mathbf{E}\log
\bigl|\lambda_j^{(\varepsilon)}-z-r\xi\bigr|I_{\Omega_n(z)}\nonumber\\
&=&
-\frac1n\sum_{j=1}^n\mathbf{E}\log\bigl(s_j\bigl(\mathbf X^{(\varepsilon
)}(z,r)\bigr)\bigr)I_{\Omega_n(z)}\\
&=&-\int_{n^{-B}}^{K_n+|z|}\log x\,d\mathbf{E}\,\overline
F_n(x,z,r),\nonumber
\end{eqnarray}
where $\overline F{}^{(\varepsilon)}_n(\cdot,z,r)$ is the distribution
function corresponding to the
restriction of the measure $\nu_n^{(\varepsilon)}(\cdot,z,r)$ to the
set $\Omega_n(z)$.
Introduce the notation
%
%
\begin{equation}
\overline U_{\mu}=-\int_{n^{-B}}^{K_n+|z|}\log x \,dF(x,z).
\end{equation}
Integrating by parts, we get
%
%
\begin{eqnarray}
\overline U{}^{(r)}_{\mu_n}-\overline U_{\mu}&=&-\int
_{n^{-B}}^{K_n+|z|}\frac{\mathbf{E}F_n^{(\varepsilon
)}(x,z,r)-F(z,r)}x\, dx\nonumber\\[-8pt]\\[-8pt]
&&{} +
C\sup_x\bigl|\mathbf{E}F_n^{(\varepsilon)}(x,z,r)-F(z,r)\bigr||{\log}(n^{B+1})|.\nonumber
\end{eqnarray}
This implies that
%
%
\begin{equation}\label{009}
\bigl|\overline U{}^{(r)}_{\mu_n}-\overline U_{\mu}\bigr|\le C\ln n\sup
_x\bigl|\mathbf{E}F_n^{(\varepsilon)}(x,z,r)-F(x,z)\bigr|.
\end{equation}
Note that, for any $r>0$, $|s_j^{(\varepsilon)}(z)-s_j^{(\varepsilon
)}(z,r)|\le r$. This implies that
%
%
\begin{equation}
\mathbf{E}F_n^{(\varepsilon)}(x-r,z)\le\mathbf{E}F_n^{(\varepsilon
)}(x,z,r)\le\mathbf{E}F_n^{(\varepsilon)}(x+r,z).
\end{equation}
Hence, we get
%
%
\begin{eqnarray}
&&\sup_x\bigl|\mathbf{E}F_n^{(\varepsilon)}(x,z,r)-F(x,z)\bigr|\nonumber\\[-8pt]\\[-8pt]
&&\qquad\le\sup
_x\bigl|\mathbf{E}F_n^{(\varepsilon)}(x,z)-F(x,z)\bigr|
+{\sup
_x}|F(x+r,z)-F(x,z)|.\nonumber
\end{eqnarray}
Since the distribution function $F(x,z)$ has a density $p(x,z)$ which
is bounded
(see Remark \ref{rem:01})
we obtain
%
%
\begin{equation}\label{08}\quad
\sup_x\bigl|\mathbf{E}F_n^{(\varepsilon)}(x,z,r)-F(x,z)\bigr|\le\sup
_x\bigl|\mathbf{E}F_n^{(\varepsilon)}(x,z)-F(x,z)\bigr|+Cr.
\end{equation}
Choose $r=\frac1{\sqrt{np_n}}$.
Inequalities (\ref{08}) and (\ref{supremum}) together imply
%
%
\begin{equation}\label{09}
\sup_x\bigl|\mathbf{E}\overline F{}^{(\varepsilon)}_n(x,z,r)-\overline
F(x,z)\bigr|\le C\biggl(\bigl(\varphi\bigl(\sqrt{np_n}\bigr)\bigr)^{-1/{18}}+\frac1{\sqrt{np_n}}\biggr).
\end{equation}
From inequalities (\ref{09}) and (\ref{009}) it follows that
\[
\bigl|\overline U{}^{(r)}_{\mu_n}-\overline U_{\mu}\bigr|\le C\biggl(\bigl(\varphi\bigl(\sqrt
{np_n}\bigr)\bigr)^{-1/{18}}+\frac1{\sqrt{np_n}}\biggr)\log(n^{B}).
\]
Note that
\[
\bigl|\overline U{}^{(r)}_{\mu_n}-U_{\mu}\bigr|\le\biggl|\int_0^{n^{-B}}\log
x\,dF(x,z)\biggr|\le Cn^{-B}|{\ln}(n^{-B})|.
\]

Let $\mathcal K=\{z\in\mathbb C \dvtx|z|\le R\}$ and let $\mathcal
K^{c}$ denote $\mathbb C\setminus\mathcal K$.
According to Lemma \ref{compact} with $q=18$, we have, for $k_1$ and
$R$ from Lemma \ref{compact},
%
%
\begin{equation}\label{rep1}
1-q_n:=\mathbf{E}\mu_n^{(r)}(\mathcal K^{c})\le\frac{k_1}{n}+\Pr\{
|\lambda_{k_1}|>R\}\le C(\varphi(np_n))^{-{1}/{18}}.
\end{equation}
Furthermore, let ${\overline{\mu}}_n^{(r)}$ and ${\widehat{\mu
}}_n^{(r)}$ be probability measures supported on the compact set
$K$ and $K^{(c)}$, respectively, such that
%
%
\begin{equation}\label{rep2}
\mathbf{E}\mu_n^{(r)}=q_n{\overline{\mu}}_n^{(r)}+(1-q_n){\widehat
{\mu}}_n^{(r)}.
\end{equation}
Introduce the logarithmic potential of the measure ${\overline{\mu}}_n^{(r)}$,
\[
U_{{\overline{\mu}}_n^{(r)}}=-{\int\log}|z-\zeta|\,d{{\overline{\mu
}}_n^{(r)}(\zeta)}.
\]
Similar to the proof of Lemma \ref{lem5.1} we show that
\[
\lim_{n\to\infty}\bigl|U_{\mu_n}^{(r)}-U_{{\overline{\mu
}}{}^{(r)}_n}\bigr|\le C\ln n (\varphi(np_n))^{-{1}/{19}}.
\]
This implies that
\[
\lim_{n\to\infty}U_{{\overline{\mu}}_n^{(r)}}(z)=U_{\mu}(z)
\]
for all $z\in\mathbb C$. According to equality (\ref{main}), $U_{\mu
}(z)$ is
equal to the potential of uniform distribution on the unit disc. This implies
that the measure $\mu$ coincides with the uniform distribution on the
unit disc.
Since the measures ${\overline{\mu}}_n^{(r)}$ are compactly
supported, Theorem 6.9 from \cite{saff}
and Corollary 2.2 from \cite{saff}
together imply that
%
%
\begin{equation}\label{rep3}
\lim_{n\to\infty}\overline{\mu}_n^{(r)}=\mu
\end{equation}
in the weak topology.
Inequality (\ref{rep1}) and relations (\ref{rep2}) and (\ref{rep2})
together imply
that
\[
\lim_{n\to\infty}\mathbf{E}\mu_n^{(r)}=\mu
\]
in the weak topology.
Finally, by Lemma \ref{sma1} we get
%
%
\begin{equation}
\lim_{n\to\infty}\mathbf{E}\mu_n=\mu
\end{equation}
in the weak topology.
Thus Theorem \ref{sparse} is proved.

\begin{appendix}\label{app}
\section*{Appendix}
In this appendix we collect some technical results.

\subsection*{The largest singular value}
Recall that $|\lambda_1^{(\varepsilon)}|\ge\cdots\ge|\lambda
_n^{(\varepsilon)}|$ denote the eigenvalues of the matrix $\mathbf
X^{(\varepsilon)}$ ordered via decreasing absolute values, and let
$s_1^{(\varepsilon)}\ge\cdots\ge
s_n^{(\varepsilon)}$ denote the singular values of the matrix $\mathbf
X^{(\varepsilon)}$.

We show the following:
\setcounter{lem}{0}
\begin{lem}\label{largeval}Under condition of Theorem \ref{thm0} for
sufficiently large $K\ge1$ we have
%
%
\setcounter{equation}{0}
\begin{equation}
\Pr\bigl\{s_1^{(\varepsilon)}\ge n\sqrt{p_n}\bigr\}\le C/{np_n}
\end{equation}
for some positive constant $C>0$.
\end{lem}
\begin{pf}Using Chebyshev's inequality, we get
%
%
\begin{equation}
\Pr\bigl\{s_1^{(\varepsilon)}\ge n\sqrt{p_n}\bigr\}\le\frac1{n^2p_n}
\mathbf{E}\operatorname{Tr} \bigl(\mathbf X^{(\varepsilon)}\bigl(\mathbf
X^{(\varepsilon)}\bigr)^* \bigr)\le
1/{(np_n)}.
\end{equation}
Thus the lemma is proved.
\end{pf}
\begin{lem}\label{compact}
Assume that $\max_{j,k}\mathbf{E}|X_{jk}|^2\varphi(X_{jk})\le C$
with $\varphi(x):=(\ln(1+|x|))^q$, $q\ge7$, and $\Delta_n:=\sup
_x|F_n^{(\varepsilon)}(x,z)-F(x,z)|$. Then there exists some absolute
positive constant $R$ such that
%
%
\begin{equation}
\Pr\bigl\{\bigl|\lambda_{k_1}^{(\varepsilon)}\bigr|>R\bigr\}\le(\varphi(np_n))^{-
({q-6})/({12q})},
\end{equation}
where $k_1:= [\Delta_n^{{(q+6)}/{(2q)}}n\ln{n} ]$.
\end{lem}
\begin{pf}
Let us introduce $k_0:= [\Delta_n^{{(q+6)}/{(2q)}}n ]$. Using
Chebyshev's inequality we obtain, for sufficiently large $R>0$,
\[
\Pr\bigl\{s_{k_0}^{(\varepsilon)}>R\bigr\}\le\frac{1-\mathbf
{E}F_n(R)}{k_0/n}\le\Delta_n^{({q-6})/({2q})}.
\]
On the other hand,
%
%
\begin{eqnarray}
\Pr\bigl\{\bigl|\lambda_{k_1}^{(\varepsilon)}\bigr|>R\bigr\}&\le&\Pr\Biggl\{\prod_{\nu
=1}^{k_1}\bigl|\lambda_{\nu}^{(\varepsilon)}\bigr|>R^{k_1}\Biggr\}\nonumber\\[-8pt]\\[-8pt]
&\le&\Pr\Biggl\{\prod_{\nu=1}^{k_1}s_{\nu}^{(\varepsilon)}>R^{k_1}\Biggr\}\le
\Pr\Biggl\{\frac1{k_1}\sum_{\nu=1}^{k_1}\ln{s_{\nu}^{(\varepsilon
)}}>\ln{R}\Biggr\}.\nonumber
\end{eqnarray}
Furthermore, for any value $R_1\ge1$, splitting into the events
$s_{k_0}^{(\varepsilon)}>R$ and
$s_{k_0}^{(\varepsilon)}\le R$, we get
\begin{eqnarray*}
&&\Pr\Biggl\{\frac1{k_1}\sum_{\nu=1}^{k_1}\ln{s_{\nu}^{(\varepsilon
)}}>\ln{R_1}\Biggr\}\\
&&\qquad\le\Pr\bigl\{s_{k_0}^{(\varepsilon)}>R\bigr\}+\Pr\biggl\{\frac
{k_0}{k_1}\ln{s_1^{(\varepsilon)}}+\ln{R}>\ln{R_1}\biggr\}\\
&&\qquad\le\Delta_n^{({q-6})/({2q})}+\Pr\biggl\{\ln{s_1^{(\varepsilon)}}>\frac
{k_1}{k_0}\ln{\frac{R_1}{R}}\biggr\}.
\end{eqnarray*}
Now choose $R_1:=R^2$. Thus, since $k_1/k_0\sim\ln n$,
\[
\Pr\bigl\{\bigl|\lambda_{k_1}^{(\varepsilon)}\bigr|>R\bigr\}\le\Delta_n^{
({q-6})/({2q})}+\Pr\bigl\{\ln{s_1^{(\varepsilon)}}>\ln{R} \ln{n}\bigr\}.
\]
Taking into account Lemma \ref{largeval} and inequality (\ref
{supremum}) we obtain
\[
\Pr\bigl\{\bigl|\lambda_{k_1}^{(\varepsilon)}\bigr|>R\bigr\}\le\Delta_n^{
({q-6})/({2q})}+\frac{C}{np_n}\le C(\varphi(np_n))^{-({q-6})/({12q})}
\]
for some positive constant $C>0$, thus proving the lemma.
\end{pf}
\begin{lem}\label{ap}Let $\varkappa=\max_{j,k}\mathbf
{E}|X_{jk}|^2\varphi(X_{jk})$. The following inequality holds:
%
%
\begin{equation}
\frac1{n\sqrt{np_n}}\sum_{j,k=1}^n\mathbf{E}\varepsilon
_{jk}|X_{jk}|\bigl(\bigl|T^{(jk)}_{k+n,j}\bigr|+\bigl|T^{(jk)}_{j,k+n}\bigr|\bigr)\le\frac
{C}{v^3\varphi(\sqrt{np_n})}.
\end{equation}
\end{lem}
\begin{pf}
Introduce the notation
%
%
\begin{equation}
B:=\frac1{n\sqrt{np_n}}\sum_{j,k=1}^n\mathbf{E}\varepsilon
{_jk}|X_{jk}|\bigl(\bigl|T^{(jk)}_{k+n,j}\bigr|+\bigl|T^{(jk)}_{j,k+n}\bigr|\bigr)
\end{equation}
and
%
%
\begin{eqnarray}
B_1&:=&\frac2{n^2p_n}\sum_{j,k=1}^n\mathbf{E}\varepsilon
_{jk}|X_{jk}|^2\bigl|R^{(jk)}_{k+n,j}\bigr|\bigl|R_{k+n,j}^{(jk)}-R_{k+n,j}\bigr|,\nonumber
\\
B_2&:=&\frac2{n^2p_n}\sum_{j,k=1}^n\mathbf{E}\varepsilon
_{jk}|X_{jk}|^2\bigl|R^{(jk)}_{k+n,k+n}\bigr|\bigl|R_{j,j}^{(jk)}-R_{j,j}\bigr|,\nonumber\\[-8pt]\\[-8pt]
B_3&:=&\frac2{n^2p_n}\sum_{j,k=1}^n\mathbf{E}\varepsilon
_{jk}|X_{jk}|^2\bigl|R^{(jk)}_{j,j}\bigr|\bigl|R_{k+n,k+n}^{(jk)}-R_{k+n,k+n}\bigr|,\nonumber
\\
B_4&:=&\frac2{n^2p_n}\sum_{j,k=1}^n\mathbf{E}\varepsilon
_{jk}|X_{jk}|^2\bigl|R^{(jk)}_{j,k+n}\bigr|\bigl|R_{j,k+n}^{(jk)}-R_{j,k+n}\bigr|.\nonumber
\end{eqnarray}
Since the function $|x|/\varphi(x)$ not decreasing, it follows from
inequality (\ref{rezolvent}) that
%
%
\begin{equation}
\bigl|R_{l,m}^{(jk)}-R_{l,m}\bigr|\le\frac1{v}I_{\{|X_{jk}|>\sqrt{np_n}\}
}+\frac1{v^2\varphi(\sqrt{np_n})}\varphi(X_{jk}).
\end{equation}
It is easy to check that
%
%
\begin{equation}
\max\{B_k, k=1,\ldots,8\}\le\frac{C\varkappa}{v^3\varphi(\sqrt{np_n})}.
\end{equation}
This implies that
%
%
\begin{equation}
B\le\frac{C\varkappa}{v^3\varphi(\sqrt{np_n})}.
\end{equation}
\upqed\end{pf}
\begin{lem}\label{ap1}
Let $\mu_n$ be the empirical spectral measure of the matrix $\mathbf
X$ and
$\nu_r$ be the uniform distribution on the disc of radius $r$. Let
$\mu_n^{(r)}$ be the empirical spectral measure of the matrix $\mathbf
X(r)= \mathbf
X-r\xi\mathbf I$, where $\xi$ is a random variable which is uniformly
distributed on the unit disc. Then the measure $\mathbf{E}\mu
_n^{(r)}$ is the
convolution of the measures $\mathbf{E}\mu_n$ and $\nu_r$, that is,
%
%
\begin{equation}
\mathbf{E}\mu_n^{(r)}=(\mathbf{E}\mu_n)*(\nu_r).
\end{equation}
\end{lem}
\begin{pf}
Let $J$ be a random variable which is uniformly distributed on the set
$\{1,\ldots,n\}$. Let $\lambda_1,\ldots,\lambda_n$ be the
eigenvalues of the
matrix $\mathbf X$. Then $\lambda_1+r\xi,\ldots,\lambda_n+r\xi$ are
eigenvalues of the matrix $\mathbf X(r)$. Let $\delta_x$ be denote the Dirac
measure. Then
%
%
\begin{equation}
\mu_n=\frac1n\sum_{j=1}^n\delta_{\lambda_j}
\end{equation}
and
%
%
\begin{equation}
\mu_n^{(r)}=\frac1n\sum_{j=1}^n\delta_{\lambda_j+r\xi}.
\end{equation}
Denote by $\mu_{nj}$ the distribution of $\lambda_j$.
Then
%
%
\begin{equation}
\mathbf{E}\mu_n=\frac1n\sum_{j=1}^n\mu_{nj}
\end{equation}
and
%
%
\begin{equation}
\mathbf{E}\mu_n^{r}=\frac1n\sum_{j=1}^n\mu_{nj}*\nu_r=
\Biggl(\frac1n\sum_{j=1}^n\mu_{nj} \Biggr)*(\nu_r)=(\mathbf{E}\mu
_n)*(\nu_r).
\end{equation}
Thus the lemma is proved.
\end{pf}

Let
%
%
\begin{equation}
f_n^{(r)}(t,v)=\int_{-\infty}^{\infty}\int_{-\infty}^{\infty}\exp
\{itx+
ivy\} \,dG_n^{(r)}(x,y)
\end{equation}
and
%
%
\begin{equation}
f_n(t,v)=\int_{-\infty}^{\infty}\int_{-\infty}^{\infty}\exp\{itx+
ivy\}\,dG_n(x,y),
\end{equation}
where
%
%
\begin{equation}
G_n^{(r)}(x,y)=\frac1n\sum_{j=1}^n\Pr\{\operatorname{Re}{\lambda
_j+r\xi}\le
x,\operatorname{Im}{\lambda_j+r\xi}\le y\}
\end{equation}
and
%
%
\begin{equation}
G_n(x,y)=\frac1n\sum_{j=1}^n\Pr\{\operatorname{Re}{\lambda_j}\le
x,\operatorname{Im}{\lambda_j}\le
y\}.
\end{equation}
Denote by $h(t,v)$ the characteristic function of the joint
distribution of
the real and imaginary parts of $\xi$,
%
%
\begin{equation}
h(t,v)=\int_{-\infty}^{\infty}\int_{-\infty}^{\infty}\exp\{iux+
ivy\}\,
dG(x,y).
\end{equation}
\begin{lem}\label{ap2}
The following relations hold
%
%
\begin{equation}
f_n^{(r)}(t,v)=f_n(t,v)h(rt,rv).
\end{equation}
If for any $t,v$ there exists $\lim_{n\to\infty}f_n(t,v)$, then
%
%
\begin{eqnarray}
\lim_{r\to0}\lim_{n\to\infty}f_n^{(r)}(t,v)&=&\lim_{n\to\infty
}\lim_{r\to0}
f_n^{(r)}(t,v)\nonumber\\[-8pt]\\[-8pt]
&=&\lim_{n\to\infty}f_n(t,v).\nonumber
\end{eqnarray}
\end{lem}
\begin{pf}The first equality follows immediately from the independence
of the random variable
$\xi$ and the
matrix $\mathbf X$. Since $\lim_{r\to0}h(rt,rv)=h(0,0)=1$ the first
equality implies the second one.
\end{pf}
\begin{lem}[(\cite{GT03}, Lemma $2.1$)]\label{ap3} Let $F$ and $G$ be
distribution functions with
Stieltjes transforms $S_F(z)$
and $S_G(z)$, respectively. Assume that
$\int_{-\infty}^{\infty}|F(x)-G(x)|\,dx<\infty$. Let $G(x)$ have a bounded
support $J$ and density bounded by some constant $K$. Let $V>v_0>0$ and
$a$ be positive numbers such that
\[
\gamma=\frac1{\pi}\int_{|u|\le a}\frac1{u^2+1} \,du > \frac34.
\]
Then there exist some
constants $C_1, C_2, C_3$ depending on $J$ and $K$ only such that
%
%
\begin{eqnarray}\qquad\quad
{\sup_x}|F(x)-G(x)|&\le& C_1 \sup_{x\in J}
\int_{-\infty}^x|S_F(u+iV)-S_G(u+iV)| \,du\nonumber\\[-8pt]\\[-8pt]
&&{} + \sup_{u\in J}
\int_{v_0}^V |S_F(u+iv)-S_G(u+iv)|\,dv+ C_3 v_0.\nonumber
\end{eqnarray}
\end{lem}
%
%
\begin{lem}\label{c01}Let $X_{jk}$, $1\le j,k\le n$, be independent
complex random variables with $\mathbf{E}X_{j,k}=0$ and $\mathbf
{E}\vert X_{j,k}\vert^2=1$. Assume furthermore that
\[
\max_{j,k}\mathbf{E}|X_{jk}|^2 I_{\{|X_{jk}|>M\}}\to0\qquad\mbox{for $M\to+\infty
$}.
\]
Then we have, for some positive $r_0$ and $\eta_0$,
\[
\sup_{u\in\mathbb C}\max_{j,k}\Pr\{\vert X_{jk}-u\vert<{\eta_0}\}
\le r_0<1.
\]
\end{lem}
\begin{pf}
First we note, that there exists a positive number $M$ such that
\[
\min_{j,k}\mathbf{E}\bigl(\vert X_{jk}\vert^2I_{\{|X_{jk}|\le M\}}\bigr)>\frac78.
\]
Let ${\eta_0}$ be a small positive number. For $|u|>M+{\eta_0}$ we have
%
%
\begin{eqnarray}\label{eq:c01:1}
\Pr\{|X_{jk}-u|\ge{\eta_0}\}&\ge&\Pr\{|X_{jk}|\le M\}\ge\frac
1{M^2}\mathbf{E}\bigl(\vert X_{jk}\vert^2I_{\{|X_{jk}|\le
M\}}\bigr)\nonumber\\[-8pt]\\[-8pt]
&>&\frac7{8M^2}.\nonumber
\end{eqnarray}
Consider now $|u|\le M+{\eta_0}$. Then
%
%
\begin{eqnarray}\label{eq:c01:2}
\Pr\{|X_{jk}-u|\ge{\eta_0}\}&\ge&\mathbf{E}\bigl(I_{\{2M+{\eta_0}\ge
|X_{jk}-u|\ge{\eta_0}\}}\bigr)\nonumber\\
&\ge&\frac1{4M^2}\mathbf
{E}\bigl(|X_{jk}-u|^2I_{\{2M+{\eta_0}\ge|X_{jk}-u|\ge{\eta_0}\}
}\bigr)\nonumber\\
&\ge&\frac1{4M^2}\bigl(1-\mathbf{E}\bigl(|X_{jk}-u|^2I_{\{|X_{jk}-u|<{\eta_0}\}
}\bigr)\nonumber\\
&&\hspace*{25.8pt}{}-\mathbf{E}\bigl(|X_{jk}-u|^2I_{\{|X_{jk}-u|>2M+{\eta_0}\}}\bigr)\bigr)\\
&\ge&\frac1{4M^2}\bigl(1-{\eta_0^2}-\mathbf{E}\bigl(|X_{jk}-u|^2I_{\{|X_{jk}|>M\}
}\bigr)\bigr)\nonumber\\
&\ge&\frac1{4M^2} \biggl(\frac34-{\eta_0^2}-\frac{|u|^2}{4M^2}
\biggr)\nonumber\\
&\ge&\frac1{16M^2} \biggl(3-4{\eta_0^2}- \biggl(1+\frac{{\eta_0^2}}{M}
\biggr)^2 \biggr).\nonumber
\end{eqnarray}
Combining inequalities (\ref{eq:c01:1}) and (\ref{eq:c01:2}) we
obtain the claim.
\end{pf}
\end{appendix}

\section*{Acknowledgments}
The authors would like to thank Terence Tao for
drawing their attention
to a gap in a previous version of the paper and Dmitry Timushev for a careful
reading of this manuscript.

%

%
\printaddresses

\end{document}